\documentclass[10pt,twocolumn]{IEEEtran}

\IEEEoverridecommandlockouts

\usepackage[utf8]{inputenc}

\usepackage{graphicx}
\usepackage{amsmath}
\usepackage{amssymb}
\usepackage{amsfonts}
\usepackage{mathrsfs}
\usepackage{dsfont}
\usepackage{color}

\usepackage[usenames,dvipsnames]{xcolor}
\usepackage{epstopdf}

\newtheorem{theorem}{Theorem}
\newtheorem{definition}{Definition}
\newtheorem{lemma}{Lemma}
\newtheorem{proposition}{Proposition}

\newtheorem{design}{Design Condition}

\newcounter{remark}
\newenvironment{remark}[1][]{\refstepcounter{remark} \par \medskip \noindent%
    \textbf{Remark~\theremark #1} \rmfamily}{\medskip}

\newcommand{\vect}[1]{\boldsymbol{#1}} 

\usepackage{enumerate}
\usepackage{balance}

\newcommand{\ak}[1]{{\color{black}{#1}}}
\newcommand{\ka}[1]{{\color{black}{#1}}}
\newcommand{\ic}[1]{{\color{black}{#1}}}
\newcommand{\icc}[1]{{\color{black}{#1}}}
\newcommand{\akk}[1]{{\color{black}{#1}}}
\newcommand{\icl}[1]{{\color{black}{#1}}}
\newcommand{\ank}[1]{{\color{black}{#1}}}
\usepackage[mathscr]{euscript}
\usepackage{enumitem}

\providecommand{\keywords}[1]{\textbf{\textit{Index terms---}} #1}
\allowdisplaybreaks

\title{
{Primary frequency regulation in power grids with on-off loads: chattering, limit cycles and {convergence to} optimality}
}

\author{Andreas Kasis\thanks{Andreas Kasis and Ioannis Lestas are with the Department of Engineering, University of Cambridge, Trumpington Street, Cambridge, CB2 1PZ, United Kingdom; e-mails: ak647@cam.ac.uk, icl20@cam.ac.uk}, Nima Monshizadeh\thanks{Nima Monshizadeh is with the Engineering and Technology Institute, University of Groningen, Nijenborgh 4, 9747AG, Groningen, The Netherlands.
email: n.monshizadeh@rug.nl} and Ioannis Lestas\thanks{A preliminary version of this study has appeared in \cite{kasis2018stability}. This manuscript  includes  additional results related to the stability and optimality properties of the considered system, further discussion and simulations and the  analytic proofs of all the main results.}
\thanks{This work was supported by ERC starting grant 679774.}
}

\begin{document}

\maketitle


\begin{abstract}
Load side participation can provide valuable support to the power network in case of urgencies.
On many occasions, loads are naturally represented by on and off states.
{{However, {{the use of {on-off} loads for frequency control can}} lead to chattering and undesirable limit cycle behavior, which are issues that need to  be resolved for {such} loads to be} used for network support.}
{This paper considers the problem of primary frequency regulation with ancillary service from on-off loads in power networks and establishes conditions that {lead to} convergence {guarantees and an appropriate power allocation within the network}.}
In particular, in order to assist existing frequency control mechanisms, we consider loads that switch when prescribed frequency thresholds are exceeded.
 {Such} control policies are prone to chattering, which limits their practicality.
 To resolve this issue, we consider loads that follow a decentralized hysteretic on-off policy, and show that chattering is not observed within such {a} setting.
 {Hysteretic loads may {exhibit}, however,} limit cycle behavior, which is undesirable.
{To {address} this, we propose {an adapted} hysteretic control scheme {for which we provide {convergence} guarantees.}}
Furthermore, we consider a mixed-integer optimization problem for {power allocation}
and {propose a suitable design of the control policy  such that the {cost incurred at equilibrium is within $\epsilon$ from the optimal cost,} 
{providing} a {non conservative value for $\epsilon$.}}
 The practicality of our analytic results is demonstrated with numerical simulations on the Northeast Power Coordinating Council (NPCC) 140-bus system.
%
%
\end{abstract}


\section{Introduction}

\textbf{Motivation and literature review:} Renewable sources of generation are expected to increase their penetration
in power networks {over} the next years \cite{lund2006large},~\cite{ipakchi2009grid}. This will result in {an} increased intermittency in the generated power, endangering power quality and potentially  the stability of the power network.
This encourages further study of the stability properties of the power grid.
Controllable loads are considered to be a way to counterbalance intermittent generation, due to their ability to provide fast response at urgencies, \ak{e.g. when there is a sudden generation/demand change or a failure in infrastructure}, by adapting their demand accordingly.
 In {recent} years, various research studies considered controllable demand as a means to {support} primary~\cite{molina2011decentralized}, \cite{kasis2016primary}, {\cite{devane2016primary}}, 
 and secondary~\cite{mallada2017optimal}, \cite{trip2016internal},  \cite{trip2018distributed}, \cite{kasis2017stability},~\cite{alghamdi2018conditions} frequency control mechanisms, with respective objectives to ensure that generation and demand are balanced and that {the} frequency converges to its nominal value (50Hz or 60Hz). Furthermore, an issue of fairness in the power allocation between controllable loads is raised if those are to be incorporated in power \ic{networks}. This problem has been pointed out in various studies \cite{kasis2017stability}, \cite{zhao2014design}, \cite{dorfler2016breaking}, \cite{kasis2020distributed}. Attempts to address this problem resulted in crafting the equilibrium of the system to coincide with the global solutions of appropriate optimization problems that ensured a fair power allocation.


On many occasions, loads are naturally represented by a discrete set of possible demand values, e.g. on and off states, and hence a continuous representation does not suffice for their study. The possible switching nature of loads has been taken into account in \cite{kasis2017secondary}, \cite{kasis2017secondary_arx}, which considered on-off loads that switched when some frequency deviation was reached in order to support the network at urgencies within the secondary frequency control timeframe. Furthermore, \cite{liu2016non} considered two switching modes of operation for loads (at nominal and urgent situations), where controllable load inputs were determined from the local deviations in frequency.
The (temperature dependent) on-off behavior of loads has been also pointed out in several studies \cite{short2007stabilization}, \cite{angeli2012stochastic}, \cite{kasis2019frequency}, \cite{aunedi2013economic}, where various control schemes for improved performance have been explored.
The study of  on-off loads with the ability to
provide support to the power network is therefore of {{major significance}} for the development of demand response schemes.
 Furthermore, the fast response required to provide ancillary support at urgencies coincides with the primary frequency control timeframe, which makes its study highly relevant for this purpose.



\textbf{Contribution:}
{This paper considers the problem of ensuring stability of the network, and optimality of the power allocation, when on-off loads contribute to primary frequency control.
{This is a problem that is significantly more involved relative to the case where only continuous generation/loads are present, since the on-off nature of the loads renders the underlying \ak{dynamical} system a hybrid system.}
Furthermore, as it will be discussed within the paper, the lack} of integral action in primary frequency control, which results to
{a non-zero steady state frequency deviation}, {further complicates the analysis} by raising problems related with the existence of equilibria and the presence of limit cycles.
{The} on-off nature of loads introduces {also} challenges in {achieving an optimal power allocation, as the corresponding network optimization problem is a mixed-integer programming problem that is NP-hard} (e.g. \cite{schrijver1998theory}).

{Our study} considers frequency dependent on-off loads {that turn on/off when
 {sufficiently large} frequency deviations occur}, within the primary frequency control timeframe, building {upon} ideas presented in\footnote{{Note that \cite{kasis2017secondary_arx} considers secondary frequency control, i.e the frequency deviation is zero at steady state and thus the loads do not contribute  \ka{to the asymptotic behavior of the system}. In this paper we consider instead primary control, which as discussed in the previous paragraph is more involved, since loads can contribute at equilibrium, which complicates the stability {and optimality analysis}. }
 }~\cite{kasis2017secondary_arx}.
We first show that the inclusion of loads that switch at a prescribed frequency does not compromise the stability of the power network, and improves {the} frequency performance.
However, such control policies {can lead to chattering, which limits} their practicality.
{A {classical} approach to resolve this is to consider hysteresis in the on-off load dynamics.
However, the coupling between frequency dynamics and load behavior  in conjunction with the discontinuous nature of the loads {can lead to cases where equilibirium points do not exist or limit {cycles occur}.}}

 {A main result of this paper is {to propose} an adapted hysteretic control scheme for on-off loads that resolves 
 {such stability issues}
 using aggregate demand measurements.}
{In particular, stability} guarantees are provided for this {scheme}, and the absence of  chattering {is also} 
analytically proven.

{A further objective of this study is to consider the problem of power allocation within the network at steady state, by requiring this to be the solution of an appropriately constructed optimization problem. Due to the discrete nature of the loads this is a mixed-integer optimization problem which is known to be NP-hard.}
{Within the paper we {propose a} 
control policy such that the {cost incurred} at equilibrium is guaranteed to be within $\epsilon$ of the optimal  {cost}, where $\epsilon$ is shown to be non-conservative.}

{A distributed mechanism for obtaining the required demand measurements is also proposed 
and we show} that the presented stability and optimality properties of the system are {unaltered with this policy.}

{Our stability and optimality analysis is numerically verified through simulations on the NPCC 140-bus system which} demonstrate that the inclusion of frequency dependent on-off loads provides improved performance and optimal steady state power allocation.

 Our {main} contribution can be summarized as follows:
 \begin{itemize}
 \item[1.] {We propose {control schemes} for on-off hysteretic loads  that {lead to  convergence guarantees, which translate to lack of limit cycles and chattering}.}
 \item[2.] {We consider a mixed-integer optimization problem {for power allocation}
 and provide design conditions for} hysteretic loads such that the {cost incurred at the resulting equilibrium points is within $\epsilon$ from the optimal cost} to this problem, providing a non-conservative {value for} $\epsilon$.
 \end{itemize}
%

\textbf{Paper structure:} The structure of the paper is as follows: Section \ref{Notation} includes some basic notation and in Section \ref{sec:Network_model} we present the power network model.
In Section \ref{Sec:switch} we consider controllable demand that switches on-off whenever certain frequency thresholds are met and present our results concerning network stability.
In Section \ref{Sec:Hysteresis}, we consider controllable loads with hysteretic {control policies}.
 In Section \ref{sec:alternative_scheme}, we propose a scheme to resolve the issue of potential limit cycle behavior from hysteretic loads and provide relevant asymptotic stability guarantees. {In Section} \ref{sec:alternative_scheme_opt} we extend our proposed scheme by {considering also the problem of optimal power allocation.}  
Numerical investigations of the results are provided  in Section \ref{Simulation}.
Finally, conclusions are drawn in Section \ref{Conclusion}. The proofs of the main results are provided in  \akk{the Appendix}. 

\vspace{-0.0mm}
\subsection{Notation}\label{Notation}
\vspace{-0.0mm}

Real and natural numbers are denoted by $\mathbb{R}$ and \ank{$\mathbb{Z}$ and} the set of n-dimensional vectors with real entries is denoted by $\mathbb{R}^n$.
The set of natural numbers including zero is denoted by $\mathbb{N}_0$
 and the {sets of positive and non-negative real numbers by $\mathbb{R}_{> 0}$ and $\mathbb{R}_{\geq 0}$.}
 Furthermore, the set of integer numbers is denoted by $\mathbb{Z}$.
We use $\vect{0}_n$ and $\vect{1}_n$  to denote $n$-dimensional vectors with all elements equal to $0$ and $1$ respectively.
The cardinality of a discrete set $\Sigma$ is denoted by $|\Sigma|$.
Moreover, we denote the collection of subsets of $\mathbb{R}^n$ by $\ka{\mathcal{P}}  (\mathbb{R}^n)$.

\vspace{-0.0mm}
\section{Network model}\label{sec:Network_model}

We describe the power network model by a connected graph $(N,E)$ where $N = \{1,2,\dots,|N|\}$ is the set of buses and $E \subseteq N \times N$ the set of transmission lines connecting the buses.
Furthermore, we use $(i,j)$ to denote the link connecting buses $i$ and $j$ and assume that the graph $(N,E)$ is directed with an arbitrary orientation, so that if $(i,j) \in E$ then $(j,i) \notin E$.
For each $j \in N$, we use
\ak{$N^i_j$ and $N^o_j$}
to denote the sets of buses that are predecessors and successors of bus $j$ respectively. It is important to note that the form of the dynamics in~\eqref{sys1}--\eqref{sys2} below is unaltered by any change in the graph ordering, and all of our results are independent of the choice of direction.
The following assumptions are made for the network:
\begin{enumerate}[label =  \arabic*)]
\item Bus voltage magnitudes are $|V_j| = 1$ \ak{per unit} for all $j \in N$.
\item Lines $(i,j) \in E$ are lossless and characterized by their susceptances $B_{ij} = B_{ji} > 0$.
\item Reactive power flows do not affect bus voltage phase angles and frequencies.
\item Relative phase angles are sufficiently small such that the approximation $\sin \eta_{ij} = \eta_{ij}$ is valid.
\end{enumerate}

Conditions $1)$ to $3)$ have been widely used in studies associated with frequency control in power networks, \ka{e.g.} \cite{mallada2017optimal}, \cite{trip2018distributed}, \cite{kasis2017secondary_arx}.
These assumptions are valid in medium to high voltages where transmission lines are dominantly inductive and voltage variations are small.
Condition $4)$ is valid when the network operates under nominal conditions, where relative phase angles are small.
Note that although the theoretical analysis relies on the above assumptions, the numerical simulations in Section \ref{Simulation}, that verify the results in this paper, make use of a full complexity model of the power network.

We use swing equations to describe the rate of change of frequency at each bus. This motivates the following system dynamics (e.g. \cite{Bergen_Vittal}),
\begin{subequations} \label{sys1}
\begin{equation}
\dot{\eta}_{ij} = \omega_i - \omega_j, \; (i,j) \in E, \label{sys1a}
\end{equation}
\begin{equation}
 M_j \dot{\omega}_j = - p_j ^L + p_j^M - (d^c_j + d^u_j)
 - \sum_{\ak{k \in N^o_j}} p_{jk} + \sum_{\ak{i \in N^i_j}} p_{ij}, j\in N, \label{sys1b}
 \end{equation}
 \begin{equation}
p_{ij}=B_{ij} \eta _{ij}, \; (i,j) \in E. \label{sys1d}
\end{equation}
\end{subequations}

In system~\eqref{sys1} the  state $\omega_j$ represents
the deviation from the nominal value\footnote{A nominal value is defined as an equilibrium of \eqref{sys1} with frequency equal to 50Hz (or 60Hz).} of {the} frequency at bus $j$.
\ak{Moreover, the state $\eta_{ij}$ represents the power angle difference\footnote{The variables $\eta_{ij}$ represent the angle difference between buses $i$ and $j$, i.e. $\eta_{ij} = \theta_i - \theta_j$, where $\theta_j$ is the angle at bus $j$. The angles themselves must also satisfy $\dot{\theta}_j = \omega_j$ at all $j \in N$. This equation is omitted in \eqref{sys1} since the power transfers are functions of the phase differences only.} and  $p_{ij}$
the power transmitted from bus $i$ to bus $j$.}
\ak{In addition,  the mechanical power injection at bus $j$ is denoted by $p^M_j$.}
The \ak{variable} $d^c_j$ represents the deviation from the nominal value\footnote{A nominal value of the controllable demand, $d^{c,nom}$, is a constant demand value selected by the users. The variable $d^c_j$ represents the deviation of the actual controllable demand from $d^{c,nom}_j$. For convenience in presentation, $d^{c,nom}_j$ is incorporated in $p^L$.}
 of controllable demand at bus $j$.
The \ak{variable} $d^u_j$  represents the uncontrollable frequency-dependent load and generation damping present at bus $j$.
The constant $M_j > 0$ denotes the generator inertia. {Moreover, the constant $p^L_j$ denotes
the frequency-independent load and the nominal value of the
controllable load at bus $j$, and $\ell =\vect{1}^T_{|N|} p^L$ its aggregate value throughout the network.}

\vspace{-0.0mm}
\subsection{Generation and uncontrollable demand dynamics}

We {consider} generation and frequency dependent uncontrollable demand dynamics described by
\begin{subequations}\label{sys2}
\begin{gather}
\tau_j \dot{p}^M_j = -(p^M_j + \alpha_j\omega_j), \; j \in N, \label{sys2_pm}
\\
d^u_j = A_j \omega_j, \; j \in N, \label{sys2_du}
\end{gather}
\end{subequations}
where $\tau_j > 0$ are time constants and $A_j > 0$ and $\alpha_j > 0,$ $j~\in~N,$ are damping and droop coefficients respectively.
Note that the analysis carried in this paper is still valid for more general generation/demand dynamics, including cases of nonlinear and higher order dynamics, provided certain input-output conditions hold, as shown in  \cite{kasis2017secondary_arx}, \cite{kasis2017primary}, \cite{monshizadeh2019secant}.
We choose to use the simple first order generation and static uncontrollable demand dynamics for simplicity and to avoid a shift in the focus of the paper from on-off loads.

\vspace{-0.0mm}
\section{On-off loads}\label{Sec:switch}

Within this section, we {consider} frequency dependent on-off loads that respond to frequency deviations by switching to an appropriate state in order to aid the network at urgencies.

\begin{figure}[t!]
\centering
\includegraphics[scale = 0.65,clip=true]{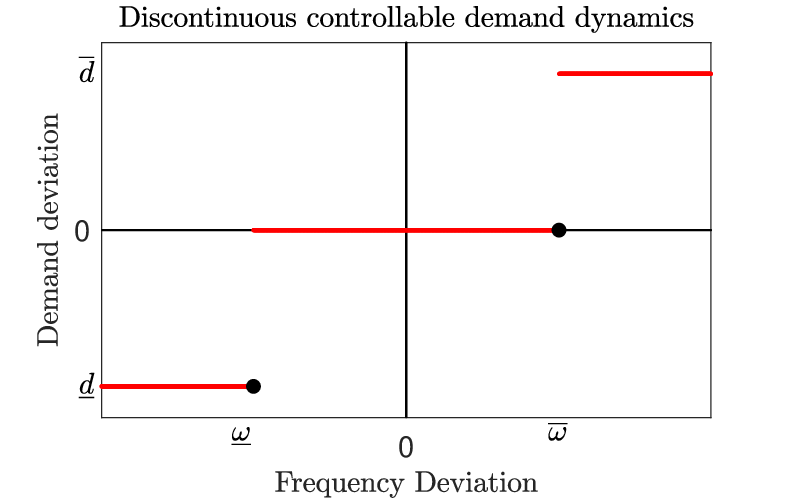}
\caption{On-off controllable demand deviations as described by \eqref{sys2_dc}.}
\label{figure_switch}
\end{figure}

The considered controllable demand dynamics are described by
\begin{equation}\label{sys2_dc}
d^c_j = f^d_j(\omega_j) =
\begin{cases}
\overline{d}_j, \quad \omega_j > \overline{\omega}_j, \\
0, \quad \hspace{4pt} \underline{\omega}_j < \omega_j \leq \overline{\omega}_j, \\
\underline{d}_j, \quad \omega_j \leq \underline{\omega}_j,
\end{cases} j \in N,
\end{equation}
where  $-\infty < \underline{d}_j \leq 0 \leq \overline{d}_j  < \infty $,   $\overline{\omega}_j > 0 > \underline{\omega}_j$ for all $j \in N$ and $f^d_j : \mathbb{R} \rightarrow \mathbb{R}$ is a discontinuous map from frequency to controllable demand at bus $j$. The static {map} in \eqref{sys2_dc} is depicted on Figure \ref{figure_switch}. Note that  \eqref{sys2_dc} may be trivially extended to include more discrete values, that would possibly respond to higher frequency deviations. The extension has been omitted for simplicity.

{To cope with the discontinuous behavior of loads and allow well defined solutions of \eqref{sys1}--\eqref{sys2_dc} for all times, a common approach is to relax \eqref{sys2_dc}} using a Filippov set valued map \cite{cortes2008discontinuous} as follows:
\begin{equation}\label{sys2_dc_set_valued}
F[d^c_j\ak{(\omega_j)}]  =  \begin{cases}
[0, \overline{d}_j],  \hspace{-1.55mm} \quad \omega_j = \overline{\omega}_j \\
[\underline{d}_j,0], \quad \hspace{-1.35mm} \omega_j = \underline{\omega}_j, \\
{\{f^d_j(\omega_j)\}, \text{ otherwise},}
\end{cases} j\in N.
\end{equation}

The states of the interconnected system \eqref{sys1}--\eqref{sys2_dc} are  denoted by $ x = (\eta, \omega, p^{M})$, where any variable without subscript represents a vector with all respective components.
For a compact representation of this system, consider the Filippov set valued map $Q: \mathbb{R}^n \rightarrow \ka{\mathcal{P}} (\mathbb{R}^n)$, where $n = |E| + 2|N|$, such that
\begin{equation}\label{sys_Filippov_representation}
\dot{x} \in Q(x)
\end{equation}
where
\begin{equation*}
\vspace{-0.0mm}
Q(x):= \begin{cases}
 \{\omega_i - \omega_j\}, \; (i,j) \in E, \\
\{\frac{1}{M_j} (- p_j ^L + p_j^M - A_j\omega_j -v_j - \sum_{\ak{k \in N^o_j}} p_{jk} \\+ \sum_{\ak{i \in N^i_j}} p_{ij} : v_j \in F[d^c_j] \} , j \in N,
 \\
 \{-\frac{1}{\tau_j}(p^M_j + \alpha_j \omega_j)\}, \; j \in N.
\end{cases}
\end{equation*}
 This representation allows the discontinuous frequency derivatives to be well-defined at all points.

{
For the analysis of system \eqref{sys1}--\eqref{sys2_dc}, we will be considering its Filippov solutions (e.g. \cite{cortes2008discontinuous}). In particular, a Filippov solution of \eqref{sys1}--\eqref{sys2_dc} on an interval $[0,t_1]$ is an absolutely continuous map $x(t)$,  $x: [0,t_1]\rightarrow\mathbb{R}^n$
that satisfies \eqref{sys_Filippov_representation} for almost all $t\in[0,t_1]$.}
Filippov solutions are often employed to analyze  {discontinuous systems},  as a means to
overcome the complications associated with the discontinuity of the vector field.

\vspace{-0.0mm}
\subsection{{Equilibrium and existence of solutions}}


{We
describe below} {what} is meant by an equilibrium of the interconnected system~\eqref{sys_Filippov_representation}.

\begin{definition} \label{eqbrdef}
The constant $x^* = (\eta^*, \omega^*, p^{M,*})$ defines an equilibrium of the system~\eqref{sys_Filippov_representation} if $\vect{0}_n \in Q(x^*)$.
\end{definition}

Note that the corresponding equilibrium value of the vector $d^{u,*}$ follows directly from $\omega^*$. Similarly the steady state controllable demand $d^{c,*}$ satisfies $d^{c,*}_j \in F[f^d_j(\omega_j^*)], j \in N$. Furthermore, note that  an equilibrium of \eqref{sys_Filippov_representation} always {exists}.


In order to study the behavior of \eqref{sys_Filippov_representation}, it is necessary to address the existence of solutions, which is stated in the following lemma, proven in \akk{the Appendix}.


\begin{lemma}\label{Uniqueness_Existence}
 There exists a Filippov solution of {system \eqref{sys1}--\eqref{sys2_dc}}  from any initial condition $x_0 =  (\eta(0), \omega(0), p^M(0))$ $\in \mathbb{R}^n$.
\end{lemma}

\vspace{-0.0mm}
\subsection{Stability analysis}

{We now present} the main result of this section, with the proof provided in \akk{the Appendix}.

\begin{theorem} \label{convthm}
The Filippov solutions of {system~\eqref{sys1}--\eqref{sys2_dc}} {converge for all initial conditions} to {an equilibrium point,} as defined in Definition~\ref{eqbrdef}.
\end{theorem}

The above theorem shows that all Filippov solutions {of~\eqref{sys1}--\eqref{sys2_dc} converge} to {an equilibrium point} of the system. It therefore demonstrates that the inclusion of controllable loads described by \eqref{sys2_dc} does not compromise the stability of the system.
{However, convergence of Filippov solutions to an equilibrium point does not rule out chattering,  as explained below, which is {a problematic behavior}. Nevertheless, Theorem \ref{convthm} provides valuable intuition on the convergence properties of the system, {used in the derivations of the results presented in the following sections}.}

\subsection{Chattering}\label{sec_chattering}

A possibility when discontinuous systems are involved, is the occurrence of infinitely many switches within some finite time, a phenomenon known as chattering (e.g. \cite{goebel2012hybrid}). Such behavior {is not} acceptable in practical implementations and should be avoided.

Chattering may occur in controllable loads, as shown in simulations in Section \ref{Simulation}.
\ak{Such behavior may occur when \ic{the component of the vector field that gives $\dot \omega_j$ for some bus $j$, changes sign  when $\omega_j$ is on either side of a point of discontinuity  $\overline{\omega}_j$ or $\underline{\omega}_j$,
such that the vector field when $\omega_j$ is on either side of this point is pointing towards this point.}}
For \ka{example}, when $0 < M_j \dot{\omega}_j < \bar{d}_j$ at some time \ic{instance} where ${\omega_j = \overline{\omega}_j}$,  then $ \dot{\omega}_j < 0$ when {a switch from off to on occurs,} which in turn \ic{causes the} frequency to decrease. This change in derivative sign will cause an infinite number of switches within some finite time, resulting {in} the aforementioned chattering behavior.

\vspace{-0mm}
\section{Hysteresis on controllable loads}\label{Sec:Hysteresis}

In this section we discuss how on-off load dynamics can be modified in order to ensure that no chattering  will occur.
To this end, we consider the use of {hysteresis} such that controllable loads switch on when a particular frequency is reached and switch off at a different frequency that is closer to the nominal one.
\ak{Such dynamics render the combined power network a hybrid system.
A formal definition of a hybrid system and its solutions are given later in this section.
The controllable load trajectories satisfy\footnote{\label{footnote_hysteresis}
\akk{It should be clarified that \eqref{sys_hysteresis} is \icl{stated as a} 
property \icl{of trajectories} $\sigma\icl{(t)}$ and $\omega\icl{(t)}$,  \icl{$t\in[t_1,t_2)\subset\mathbb{R}$} but \icl{does not define how these are generated.}
An exact definition of solutions, using a hybrid systems formalism, is provided later in Definition \ref{dfn_hybrid_domain_solution}.
}}
 the following property}
\begin{equation}\label{sys_hysteresis}
d^c_j\icl{(t)} = \overline{d}_j \sigma_{j}\icl{(t)}, \quad
{\akk{\sigma_{j}^+(t)} \in \begin{cases}
\{1\}, \qquad \hspace{1.5mm}  \omega_j\icl{(t)} >  \omega^1_j\\[1mm]
\{0\},  \qquad  \hspace{1.5mm}  \omega_j\icl{(t)} <  \omega^0_j\\[1mm]
\{\sigma_j(t)\}, \hspace{3mm}  \omega^0_j < \omega_j\icl{(t)} < \omega^1_j \\
 \{0, \sigma_j(t)\},  \omega_j\icl{(t)} = \omega^0_j \\
\{\sigma_j(t), 1\},  \omega_j\icl{(t)} = \omega^1_j
\end{cases}}
\end{equation}
where  $j\in N$, $\akk{\sigma_{j}^+(t) = \lim_{\epsilon \rightarrow \ank{0^+}} \sigma(t + \epsilon)}$, $\overline{d}_j > 0$ \ank{and} the frequency thresholds $\omega_j^0, \omega_j^1$, satisfy
$\omega_j^1 >\omega_j^0>0$.
\ank{Signal}
$\sigma_j\icl{(t)} \in P = \{0,1\}$ denotes the switching state for loads in bus $j \in N$
\ank{that is continuous at all times $t$ where $\omega_j(t)\notin \{\omega^0_j, \omega^1_j\}$. Also if $\sigma_j(t)$ is discontinuous at $t$, with $\omega_j(t)=\omega^0_j$ (similarly $\omega_j(t)=\omega^1_j$) then $\sigma_j^+(t)=0$ (similarly  $\sigma_j^+(t)=1)$.}
  {For generality, the control scheme \eqref{sys_hysteresis} considers two possibilities when frequency thresholds $\omega^0_j$ and $\omega^1_j$ are reached, corresponding to a switch when the frequency reaches or exceeds a particular threshold. This approach is used throughout the rest of the paper and is consistent with the widely used framework in \cite{goebel2012hybrid} for the analysis of hybrid systems.
Note that the results in Sections \ref{Sec:Hysteresis}--\ref{sec:alternative_scheme_opt} concerning convergence of solutions and absence of chattering are about all solutions of the resulting hybrid systems.
}

\ak{
\begin{remark}\label{rem_delay}
Chattering behaviour could be avoided by implementing some time delay in the switch of on-off loads, i.e. enabling a switch only when the frequency is below a particular threshold for some given time \ka{duration}.
However, it can be shown that  such schemes may result in absence of equilibria to the power system and hence in undesirable behaviors.
\akk{In addition, the imposed time delay may reduce the \icl{effectiveness} of the  ancillary services provided  from on-off loads.
}
\end{remark}
}

\begin{figure}[t]
\centering
\includegraphics[trim = 0mm 0mm 0mm 0mm, scale = 0.65,clip=true]{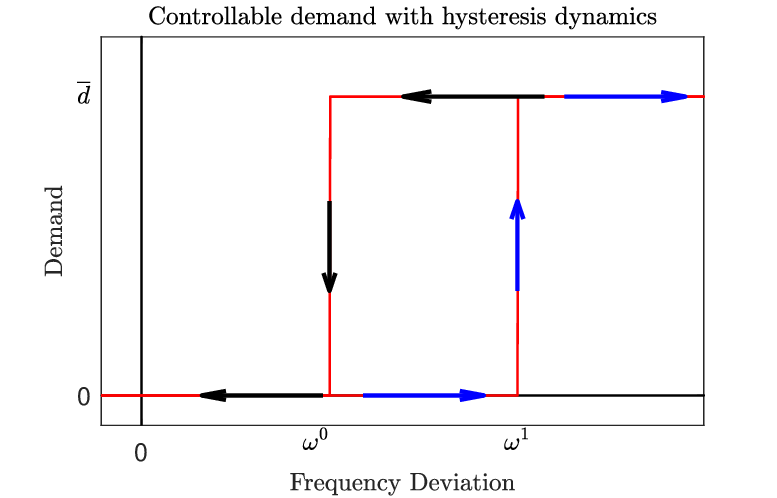}
\vspace{-1mm}
\caption{Hysteresis dynamics for controllable loads described by \eqref{sys_hysteresis}.}
\label{hysteresis_figure}
\vspace{-0mm}
\end{figure}

The dynamics in \eqref{sys_hysteresis} describe loads that switch on from off. Note that the conjugate case of loads switching off from on can also be incorporated by reversing the signs of frequency thresholds and controllable demand deviations and that all the analytic results of this paper can be trivially extended to include this case. However, we consider only loads that switch from off to on for simplicity in presentation.
 The dynamics described in \eqref{sys_hysteresis}
  can be visualized in Figure \ref{hysteresis_figure}.
 Moreover, we use {$t_{i,j}, i \in \mathbb{N}, j \in N$} to denote the time-instants where the value of $\sigma_j$ changes. Within the rest of the paper we shall adopt the notation \akk{$a^+(t) = \lim_{\epsilon \rightarrow \ank{0^+}} a(t + \epsilon)$} for any real vector $a(t)$.
\ak{For convenience in \ic{the} notation, we will refer to $a^+(t)$ by simply $a^+$.}

The behavior of system \eqref{sys1},~\eqref{sys2},~\eqref{sys_hysteresis} can be described by the states $ \zeta = (x,\sigma$), where $x = (\eta, \omega, p^{M}) \in \mathbb{R}^n$, $n = |E| + 2|N|$, is the continuous state, and $\sigma \in P^{|N|}$ the discrete state. Moreover, let $\Lambda =  \mathbb{R}^n \times P^{|N|}$ be the space where the system's states evolve.
The continuous dynamics of the system \eqref{sys1},~\eqref{sys2},~\eqref{sys_hysteresis} are described by
\begin{subequations}\label{sys4_hysteresis}
\begin{gather}
\dot{\eta}_{ij} = \omega_i - \omega_j, \; (i,j) \in E, \label{sys4a} \\
\hspace{-5mm} M_j \dot{\omega}_j = - p_j ^L + p_j^M - (\overline{d}_j\sigma_j + A_j \omega_j) \nonumber \\
\hspace{15mm} - \sum_{\ak{k \in N^o_j}} p_{jk} + \sum_{\ak{i \in N^i_j}} p_{ij}, \; j\in N, \label{sys4b} \\
 p_{ij}=B_{ij} \eta _{ij}, \; (i,j) \in E, \label{sys4c} \\
 \tau_j \dot{p}^M_j = -(p^M_j + \alpha_j\omega_j), \; j \in N, \label{sys4d} \\
 \dot{\sigma}_j = 0, j \in N, \label{sys4e}
\end{gather}
\end{subequations}
%
which is valid when $\zeta$ belongs to the set $C$ described below,
\begin{equation}\label{e:C}
C=\{ \zeta \in \Lambda:  \sigma_j \in \mathcal{I}_j(\omega_j), \;\forall j\in N\},
\end{equation}
where
\[
\mathcal{I}_j(\omega_j)= \begin{cases}
\{1\}, \qquad {\omega_j >  \omega^1_j},\\[1mm]
\{0\},  \qquad {\omega_j <  \omega^0_j},\\[1mm]
\{0, 1\}, \quad \omega^0_j \leq \omega_j \leq \omega^1_j.
\end{cases}
\]

Alternatively, when $\zeta$ belongs to the set $D = {\Lambda \setminus C \cup \underline{D}}$
where
$\underline{D} = \{ \zeta \in \Lambda:
\sigma_j \in \mathcal{I}^D_j(\omega_j), \;\forall j\in N\}
$,
and
\[
\mathcal{I}^D_j(\omega_j)= \begin{cases}
\{0\}, \qquad \omega_j =  \omega^1_j,\\[1mm]
\{1 \}, \qquad \omega_j =  \omega^0_j,
\end{cases}
\]
its components follow the discrete update depicted below
\begin{align}\label{sys4_g}
x^+ &= x, &
\akk{\sigma^+_{j}} = \begin{cases}
1, \quad \omega_j \geq  \omega^1_j,\\[1mm]
0,  \quad \omega_j \leq  \omega^0_j.
\end{cases}
\end{align}
 We can now provide the following compact representation for the hybrid system \eqref{sys1},~\eqref{sys2},~\eqref{sys_hysteresis},
\begin{subequations}\label{sys4}
\begin{gather}
\dot{\zeta} = f(\zeta), \zeta \in C, \\
\zeta^+ = g(\zeta), \zeta \in D,
\end{gather}
\end{subequations}
where $f(\zeta): C \rightarrow {C}$ and $g(\zeta): D \rightarrow {C \setminus D}$  are described by \eqref{sys4_hysteresis} and \eqref{sys4_g} respectively.  Note that $\zeta^+ = g(\zeta)$ represents a discrete dynamical system where $\zeta^+$ indicates that the next value of the state $\zeta$ is given as a function of its current value through $g(\zeta)$. Moreover, note that $C \cup D = \Lambda$.

%

\subsection{Analysis of equilibria and solutions}

In this subsection, we define and study the equilibria and solutions of \eqref{sys4}. We provide  sufficient design conditions  for the existence of equilibria of \eqref{sys4} and show that {chattering does not occur when} hysteretic dynamics {are used}.

Below, we provide {the} definition of an equilibrium of {the} system described by \eqref{sys4}.
\begin{definition}\label{eqbr_dfn_hybrid}
A point $\zeta^*$ is an equilibrium of the system described by \eqref{sys4} if it satisfies $f(\zeta^*)=0, \zeta^* \in C$ or $\zeta^* = g(\zeta^*), \zeta^* \in D$.
\end{definition}

{It should be noted that, when hysteretic loads are introduced, the system is not guaranteed to have equilibria, and hence additional conditions are required.}
%
The following theorem, proven in \akk{the Appendix}, provides {a sufficient condition under which an equilibrium to \eqref{sys4} exists}.

For the rest of the manuscript we define $\mathcal{D} = \sum_{j \in N} (\alpha_j + A_j)$.

\begin{theorem}\label{eqbr_hyst_existence}
An equilibrium point $\zeta^*$ of \eqref{sys4} exists for any $p^L$ {if} $\omega^1_j - \omega^0_j \geq \overline{d}_j / \mathcal{D}$  holds for all $j \in N$.
\end{theorem}

Theorem \ref{eqbr_hyst_existence} provides a {sufficient} {design} condition on the hysteretic dynamics which ensures that equilibria will exist for any load profile. Potential lack of equilibria results in undesirable behaviors such as limit cycles. Stability-wise, the conditions for existence of equilibria can be seen as necessary conditions for convergence to a fixed point. {Furthermore, there exist configurations where it can be shown that the condition in Theorem \ref{eqbr_hyst_existence} is also necessary, e.g. when the hysteresis region in at least one load is non-overlapping with the respective hysteresis regions of all other loads.} The physical interpretation of Theorem \ref{eqbr_hyst_existence} is that the hysteresis region of each on-off load should be no smaller than  the frequency deviation caused by its switch, which can be shown to be $\overline{d}_j / \mathcal{D}$.
\ak{Note also that  Theorem \ref{eqbr_hyst_existence} trivially holds when $\mathcal{D}$ is replaced by a known lower bound, which offers robustness to measurement uncertainty.}


Below, we provide a definition of a hybrid time domain, hybrid solution and complete and maximal solutions for systems described by \eqref{sys4}.
\ak{Note that we use the definition of a hybrid system from \cite[Dfn. 2.2]{goebel2012hybrid}.}

\begin{definition}{(\cite{goebel2012hybrid})}\label{dfn_hybrid_domain_solution}
A subset of $\mathbb{R}_{\geq 0} \times \mathbb{N}_0$ is a hybrid time domain if it is a union of a finite or infinite sequence of intervals  $[t_l, t_{l+1}] \times \{l\}$, with the last interval (if existent) possibly of the form $[t_l, t_{l+1}] \times \{l\}$, $[t_l, t_{l+1}) \times \{l\}$, or $[t_l, \infty) \times \{l\}$.
Consider a function $\zeta(t,l): K \rightarrow \mathbb{R}^n$ defined on a hybrid time domain $K$ such that for every fixed $l \in \mathbb{N}$, $t\rightarrow \zeta(t,l)$ is locally absolutely continuous on the interval $T_l = \{t: (t,l) \in K\}$.
 The function $\zeta(t,l)$ is a solution to the hybrid system $\mathcal{H} = (C,f,D,g)$ if {$\zeta(0,0) \in {{C} \cup D}$}, and for all $l \in \mathbb{N}$ such that $T_l$ has non-empty interior (denoted by ${\rm int}T_l$)
\begin{align*}
& \zeta(t,l) \in C, \text{ for all t} \in {{\rm int} T_l}, \\
& \dot{\zeta}(t,l) \in f(\zeta(t,l)), \text{ for almost all } t \in T_l,\\
& \hspace{-12mm} \text{and for all } (t,l) \in K \text{ such that } (t,l+1) \in K,\\
& \zeta(t,l) \in D, \; \zeta(t,l+1) \in g(\zeta(t,l)).
\end{align*}
A solution $\zeta(t,l)$ is complete if $K$ is unbounded. A solution $\zeta$ is maximal if there does not exist another solution $\tilde \zeta$ with time domain $\tilde K$ such that $K$ is a proper subset of $\tilde K$ and $\zeta(t,l) = \tilde{\zeta}(t,l)$ for all $(t,l) \in K$.
\end{definition}

{For convenience {in the presentation the term solutions within the paper will refer to maximal solutions\footnote{{We will also occasionally use explicitly the term maximal solutions to remind the reader of this property in cases this is technically of significance.}}.}}
The {following proposition demonstrates the} existence of solutions to \eqref{sys4} as well as of a finite dwell time between switches of states $\sigma_j$ within any compact set. {Furthermore, it establishes that all maximal solutions to \eqref{sys4} are complete.}
\begin{proposition}\label{dwell_time_lemma}
For any initial condition $\zeta(0,0) \in {\Lambda}$ there exists a complete solution to \eqref{sys4}. {All maximal solutions to \eqref{sys4} are complete.}
 Furthermore, for any complete bounded solution to \eqref{sys4}, there exists $\tau_j > 0$  such that $\min_{i \geq 1} (t_{i+1, j} - t_{i, j}) \geq \tau_j$ for any $j \in N$.
\end{proposition}

\begin{remark}
The importance of Proposition \ref{dwell_time_lemma} is that it shows that no chattering will occur for any complete bounded solution of system \eqref{sys4}. This is because for any finite time interval $\tau = \min_j \tau_j, j \in N$, the vector $\sigma$ changes at most $|N|$ times. This shows the practical advantage of \eqref{sys4} when compared to \eqref{sys_Filippov_representation}.
\end{remark}

\vspace{-0.0mm}
\subsection{Limit cycle behavior}\label{sec_limit_cycle}

Numerical simulations \ak{in Section \ref{Simulation} (see Fig. \ref{Limit_cycle_case})} demonstrate that limit cycle behavior can {occur} when the considered hysteretic loads are introduced in the network.  This is a consequence of the load on-off behavior which results to discontinuous changes in the vector field which in turn cause further switches.
{Note that
 the existence of equilibrium points does not ensure the absence of limit cycles.
 }
 In the following section we present an approach to resolve this issue.

\section{An {adapted} scheme for hysteretic loads}\label{sec:alternative_scheme}

\begin{figure}[t!]
\centering
\includegraphics[trim = 0mm 0mm 0mm 0mm, scale = 0.7,clip=true]{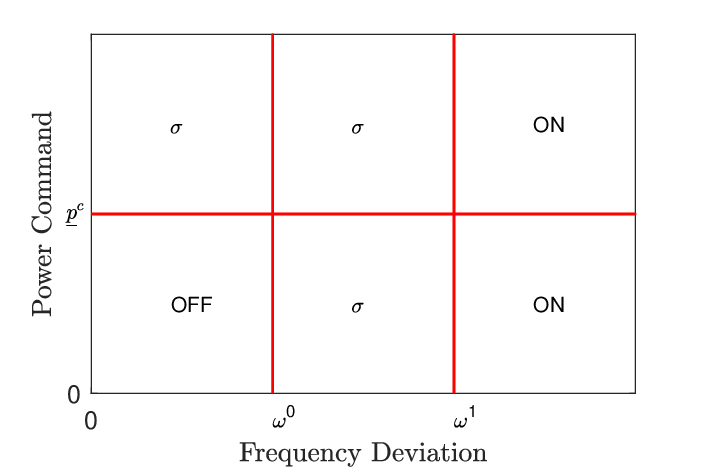}
\vspace{-0mm}
\caption{{Adapted hysteresis} scheme for controllable loads described by \eqref{sys_hysteresis_new_simple}.}
\label{hysteresis_new_figure_simple}
\vspace{-0mm}
\end{figure}

{
In this section, we discuss how  hysteretic on-off load dynamics may be modified to guarantee convergence, ruling out limit cycle behavior.
In particular, we propose a control scheme that allows two modes of operation for on-off loads;
 one that implements \eqref{sys_hysteresis}, and  a second one that allows loads to switch on  {when significant} frequency deviations are observed, providing support to the power network,  but prohibits further switches.
 The latter is in line with existing  load shedding practices where loads are switched at urgencies (e.g. \cite[Ch. 9]{machowski2011power}).
{The mode of operation of the loads is determined from the aggregate demand.
{In particular,} load shedding is implemented on an increasing portion of {on-off loads as the total demand increases}.}
  In this section, we explain how such scheme should be designed {such that on-off loads provide ancillary service to the power network without compromising its stability properties.}}
 In particular, \ak{controllable demand trajectories satisfy the following property}
\begin{equation}\label{sys_hysteresis_new_simple}
d^c_j = \overline{d}_j \sigma_{j},\hspace{-0.0mm}
\hspace{1mm}
{\akk{\sigma_{j}^+}\icl{(t)} \hspace{-0.0mm} \in  \hspace{-0.0mm}\begin{cases}
\{1\} , \qquad \hspace{2mm} \omega_j >  \omega^1_j,  \\[1mm]
\{0\},  \qquad \hspace{2mm} \omega_j <  \omega^0_j \text{ and } p^c < \underline{p}^c_j,\\[1mm]
\{\sigma_j(t)\}, \hspace{1mm}  \hspace{2mm} \begin{cases}
\omega_j^0 <\hspace{-0.0mm}   \omega_j <\hspace{-0.0mm}  \omega_j^1,   \\
\omega_j \hspace{-0.0mm} < \hspace{-0.0mm} \omega_j^0 \text{ and }   \underline{p}^c_j  <\hspace{-0.0mm}  p^c,
\end{cases} \\
\{ 0,  \sigma_j(t)\}, \begin{cases}
\omega_j =  \omega^0_j \text{ and }   p^c \leq \underline{p}^c_j,   \\
\omega_j \leq \omega_j^0 \text{ and }   p^c = \underline{p}^c_j,
\end{cases} \\
\{\sigma_j(t), 1\},  \quad \hspace{-0mm} \omega_j = \omega^1_j,
\end{cases}}
\end{equation}
where $j\in N$, $\underline{p}^c_j$ are variables available for design (see Section \ref{sec: Controller Design_simple} below), and $\overline{d}_j, \omega^0_j$ and $\omega^1_j$ are as in \eqref{sys_hysteresis}. The scheme in \eqref{sys_hysteresis_new_simple} is depicted in Figure \ref{hysteresis_new_figure_simple}.
Furthermore, $p^c$ is a power command variable given by
\begin{equation}\label{sys_power_command}
p^c  = -\ell = - \sum_{j \in N} p^L_j.
\end{equation}

\begin{remark}\label{rem:hyst}
The scheme presented in \eqref{sys_hysteresis_new_simple} {uses the power command signal {\eqref{sys_power_command}} to determine the dynamic behavior of each load}. In particular, when the power command value is above the local respective threshold $\underline{p}^c_j$, then switching from on to off is prohibited, {although loads can still switch once from off to on to support the network}.
Alternatively, when $p^c \leq \underline{p}^c_j$, then \eqref{sys_hysteresis_new_simple} reduces to \eqref{sys_hysteresis} and convergence depends on the choice of the thresholds in \eqref{sys_hysteresis_new_simple}, which are available for design.
In Section \ref{sec: Controller Design_simple} we discuss how these {thresholds should be selected} such {that convergence can be deduced.}
\end{remark}

\begin{remark}\label{rem:distributed}
\ak{An application scenario is to implement the scheme in \eqref{sys_hysteresis_new_simple}--\eqref{sys_power_command} on an aggregation of small loads at a given \icc{bus.}}
Note also that the scheme in \eqref{sys_power_command} requires knowledge of the total demand of the system \icc{which could be estimated by SCADA systems,  e.g. \cite{gomez2004power}.}
{Furthermore, note that all convergence properties presented below are retained \ak{when the magnitude of $\ell$} is replaced with a known upper bound, and hence \ak{its precise value} is not necessary for stability (see also Remark \ref{remark_exact_ell}).}
\ka{Such upper bound may be obtained using historical demand data,
or \icc{updated via signals} from the operator at slower timescales.}
\ak{In addition, the requirement for frequency measurements can be fulfilled with low cost at the load level.}
It should further be noted that the requirement for a centrally implemented controller to transmit the total demand in \eqref{sys_hysteresis_new_simple}, \eqref{sys_power_command}  is relaxed in \akk{\cite[Appendix B]{kasis2019primary}}, where we present a distributed scheme to evaluate the aggregate demand without compromising the convergence properties of the system.
\end{remark}

\subsection{Controller design}\label{sec: Controller Design_simple}

In this section we propose a way to design power command and frequency thresholds
{{such that} loads that satisfy $p^c \leq \underline{p}^c_j$ are off at steady state which, {as shown below}, allows to deduce convergence {to the set of equilibrium points}.}
The  condition concerns the power command and lower frequency thresholds $\underline{p}^c_j$ and $\omega^0_j$.
We remind that $\mathcal{D} = \sum_{j \in N} (\alpha_j + A_j)$.

\begin{design}\label{assum_pc}
The values of $\underline{p}^c_j$ and $\omega^0_j$ are chosen such that
$\underline{p}^c_j \leq \mathcal{D} \omega^0_j$ holds.
\end{design}

Design condition \ref{assum_pc} rules out the occurrence of limit cycles, as {follows from} 
Theorem \ref{conv_thm_hysteresis_simple} {in Section \ref{sec:hybrid stability}}. The scheme \eqref{sys_hysteresis_new_simple}--\eqref{sys_power_command} ensures that each load will satisfy either $p^{c} > \underline{p}^c_j$, which prohibits switching from on to off as explained in Remark \ref{rem:hyst}, or $p^{c} \leq \underline{p}^c_j$.
When the latter {occurs, switching depends on the frequency only as follows from \eqref{sys_hysteresis_new_simple}, and} Design condition  \ref{assum_pc}  guarantees that the equilibrium frequency is less than the corresponding frequency thresholds $\omega_j^0$,  {a property that is key to provide stability guarantees, as shown in the proof of Theorem~\ref{conv_thm_hysteresis_simple}.}
 The condition follows by noting that the power command and the equilibrium values of frequency depend directly on $\ell$, as shown by \eqref{sys_power_command} and \eqref{eqlb_freq} below.
\begin{equation}
\omega^* = \frac{ -\ell - \overline{d}^T \sigma^*}{\mathcal{D}}
\label{eqlb_freq}
\end{equation}
{From \eqref{eqlb_freq}, it follows that the value of $\ell$ allows to obtain an upper bound of the equilibrium frequency, attained when $\sigma^* = 0$.
Hence, Design condition \ref{assum_pc} guarantees that when $p^c < \underline{p}^c_j$, then $\omega^* < \omega^0_j$, noting that $\mathcal{D} = \frac{p^{c}}{\omega^*}\bigr|_{\sigma^* = 0}$.}
 The condition can be easily fulfilled since both $\omega_j^0$ and $\underline{p}^c_j$ are design variables.
It should  be further noted that Design condition \ref{assum_pc} requires knowledge of the aggregate droop and damping coefficients from all buses across the network. However, for the purpose of the analysis, it is sufficient to have a lower bound of $\mathcal{D}$, which offers robustness to model uncertainty.
%

{
The practical significance and non-conservativeness of the proposed scheme is demonstrated with realistic simulations in Section \ref{Simulation}, where {significant improvement in the frequency response is observed.}
}

{
\begin{remark}
An alternative approach to avoid limit cycles would be to choose the set of loads satisfying $p^c > \underline{p}^c_j$ and assign to them an arbitrary switching condition. However, such scheme would not  respond to local frequency deviations   and hence not provide an efficient ancillary service to the power network, i.e. loads could switch at buses far from a disturbance.
Furthermore, such a scheme would require central knowledge of {the  power command thresholds of all loads} and could result to increased user disutility, causing unnecessary load switch.
\end{remark}
}

\subsection{Hybrid system description}

The states {$\zeta = (x, \sigma) \in \Lambda$ can describe the behavior of the system \eqref{sys1},~\eqref{sys2},~\eqref{sys_hysteresis_new_simple},~\eqref{sys_power_command}.
%
Its continuous dynamics are described by \eqref{sys4_hysteresis} and \eqref{sys_power_command}
 when $\zeta$ belongs to the set $F$ defined below,
\begin{equation}\label{f:C_simple}
F=\{ {\zeta} \in \Lambda:  \sigma_j \in \overline{\mathcal{J}}_j(\omega_j, p^c), \;\forall j\in N\},
\end{equation}
where
\begin{equation}\label{J_set_flow_simple}
\overline{\mathcal{J}}_j(\omega_j, p^c) \hspace{-0.5mm}= \hspace{-0.5mm}\begin{cases}
\{1\}, \quad {\omega_j >  \omega^1_j,} \\[1mm]
\{0\},  \quad {
\omega_j < \omega_j^0 \text{ and }   p^c < \underline{p}^c_j, }
\\[1mm]
\{0, 1\},  \hspace{-0.5mm}\begin{cases}
\omega_j^0 \leq  \omega_j \leq \omega_j^1,  \\
\omega_j \leq \omega_j^0 \text{ and }   \underline{p}^c_j \leq p^c.
\end{cases}
\end{cases}
\end{equation}

Furthermore, when ${\zeta} \in G ={ \Lambda \setminus F \cup \tilde{G}}$,
where
$
\tilde{G} = \{ {\zeta} \in \Lambda:
\sigma_j \in \overline{\mathcal{I}}^D_j(\omega_j, p^c), \;\forall j\in N\},
$
and
\begin{equation}\label{jump_set_simple}
\overline{\mathcal{I}}^D_j(\omega_j, p^c)= \begin{cases}
\{0\}, \quad {\omega_j =  \omega^1_j},\\[1mm]
\{1 \},  \begin{cases}  \omega_j \leq \omega_j^0 \text{ and }  p^c = \underline{p}^c_j, \\
 \omega_j = \omega_j^0 \text{ and } p^c \leq \underline{p}^c_j,
\end{cases}
\end{cases}
\end{equation}
then its components follow a discrete update given by
\begin{align}\label{sys5_g_simple}
{x^+} &= x, &
\akk{\sigma_{j}^+} = \begin{cases}
1, \quad \omega_j\geq  \omega^1_j,\\[1mm]
0,  \quad \omega_j\leq  \omega^0_j \text{ and } p^c \leq \underline{p}^c_j.
\end{cases}
\end{align}
 Hence, the hybrid system \eqref{sys1},~\eqref{sys2},~\eqref{sys_hysteresis_new_simple},~\eqref{sys_power_command} can be represented by
\begin{subequations}\label{sys5_simple}
\begin{gather}
\dot{\zeta} = \tilde{f}(\zeta), \zeta \in F,
\\
\zeta^+ = \tilde{g}(\zeta), \zeta \in G,
\end{gather}
\end{subequations}
where $\tilde{f}(\zeta): F \rightarrow {F}$ and $\tilde{g}(\zeta): G \rightarrow {F \setminus G}$  follow from \eqref{sys4_hysteresis} and \eqref{sys_power_command}, and \eqref{sys5_g_simple} respectively.


\subsection{Equilibrium and solutions analysis}

Below we provide {the definition of an equilibrium of} \eqref{sys5_simple}.
\begin{definition}\label{eqbr_dfn_hybrid_2_simple}
A point $\zeta^*$ is an equilibrium of the system described by \eqref{sys5_simple} if it satisfies ${\tilde{f}}(\zeta^*)=0, \zeta^* \in F$ or $\zeta^* ={\tilde{g}}(\zeta^*), \zeta^* \in G$.
\end{definition}


 {The following proposition, proven in \akk{the Appendix}, states that equilibria of  \eqref{sys5_simple} exist} when Design condition \ref{assum_pc} {holds}.
\begin{proposition}\label{eqlbr_hysteresis_simple}
Consider the system described by \eqref{sys5_simple} and let Design condition \ref{assum_pc}
hold.
Then, an equilibrium point exists and satisfies $\zeta^* \in F$.
\end{proposition}

{{Proposition} \ref{Proposition_existence_no_Zeno_simple} below shows the existence of solutions to \eqref{sys5_simple} and of a minimum time between consecutive switches. The latter implies that no chattering {occurs}.}

\begin{proposition}\label{Proposition_existence_no_Zeno_simple}
For any initial condition {$\zeta(0,0) \in \Lambda$} there exists a complete solution to \eqref{sys5_simple}. {All maximal solutions to \eqref{sys5_simple} are complete.}
Furthermore, for any complete bounded solution to \eqref{sys5_simple}, there exists $\tau > 0$  such that $\min_{i \geq 1} (t_{i+1, j} - t_{i, j}) \geq \tau, j \in N$.
\end{proposition}


\subsection{Stability of hybrid system}\label{sec:hybrid stability}

In this section, we provide our main convergence result about system \eqref{sys5_simple}, with the proof provided in \akk{the Appendix}.

\begin{theorem}\label{conv_thm_hysteresis_simple}
Let Design condition \ref{assum_pc} hold. Then, {for all initial conditions,}
the solutions of \eqref{sys5_simple} {are bounded and converge} to {a subset} of its equilibria.
\end{theorem}


{Theorem \ref{conv_thm_hysteresis_simple} and Proposition \ref{Proposition_existence_no_Zeno_simple} show} that the inclusion of loads with dynamics described by \eqref{sys_hysteresis_new_simple} does not compromise the stability of the system, {when  Design condition \ref{assum_pc}  holds}, and neither exhibits any {chattering 
behavior.}


\section{Optimal power allocation with hysteretic loads}\label{sec:alternative_scheme_opt}

%
%
\subsection{Optimal supply and hybrid load control problem}\label{sec_optimal_supply}

{
We investigate in this section} how to adjust the generation and hybrid controllable demand to meet the step change in $p^L$  and simultaneously minimize the total cost that accounts for the extra power generated and the \ak{\ic{cost incurred  when on-off loads alter their demand}}.

Let $C_{h,j}(d^c_j)$ denote the costs incurred from deviations  $d^c_j$  in controllable demand.
The discrete nature of controllable loads suggests the following structure for the cost functions,
\begin{equation}\label{cost_function_load}
C_{h,j}(d^c_j) = \begin{cases}
0,  \;\; d^c_j = 0, \\
c^d_j, d^c_j = \bar{d}_j,
\end{cases}  j \in N,
\end{equation}
where $c^d_j > 0, j \in N$. Furthermore, we let $\frac{c_j}{2}(p_{j}^M)^2$ and $\frac{1}{2A_j}(d^{u}_j)^2$ be the costs incurred for  generation $p^M_j$
and the  change in frequency, which alters frequency dependent uncontrollable demand $d^u_j$.
The total cost is the sum of all the above costs. The problem, called the  optimal supply and hybrid load control problem (H-OSLC), is to choose the vectors $p^M$, $d^c$ and $d^u$ such that this total cost is minimized when simultaneously power balance is achieved.
\begin{equation}
\begin{aligned}
&\hspace{2em}\underline{\text{H - OSLC:}} \\
&\min_{p^M,d^c,d^u} \sum\limits_{j\in N} \Big( \frac{c_j}{2}(p_{j}^M)^2 +   C_{h,j} (d^c_{j}) + \frac{1}{2A_j}(d^{u}_j)^2  \Big) \\
&\text{subject to } \sum\limits_{j\in N} (p_j^M - d^u_j - p_j^L) = \sum\limits_{j\in  N} d^c_j, \\
& d^c_j   \in \{0,\bar{d}_j\}, j \in N.
 \label{Problem_To_Min_H}
\end{aligned}
\end{equation}
{The first constraint in~\eqref{Problem_To_Min_H} is associated with the balance between generation and demand, which is a property that needs to be satisfied at equilibrium.}
 The second constraint reflects the fact that controllable loads take discrete values, making \eqref{Problem_To_Min_H} a mixed-integer optimization problem.

\subsection{Controller design for convergence and optimality}\label{sec: Controller Design_opt}

{In this section, we propose a control scheme that allows on-off loads to provide ancillary services to the power network and simultaneously  ensures that {the cost incurred at   equilibrium  is close to the optimal cost} of \eqref{Problem_To_Min_H}.
Since the solution to \eqref{Problem_To_Min_H} determines whether a load is on or off at steady state for given aggregate demand value $\ell$, it follows that the control policy needs to allow load equilibrium values to be determined from $\ell$.
In particular, we consider two main  modes of operation for on-off loads; one where loads stay switched on at all times and a second one where loads  implement \eqref{sys_hysteresis} to provide transient support to the network, but are designed to be switched off at equilibrium.
In addition, to avoid possible chattering in the presence of noise in  demand measurements, we implement a third mode of operation which allows loads to switch once, when appreciable frequency deviations are present, but prohibits further switches, similar to \eqref{sys_hysteresis_new_simple}.
%
We then explain how appropriate selection of the threshold values
{results to a power allocation that is close to optimal.}}

 In particular, \ak{controllable demand trajectories satisfy}
\begin{equation}\label{sys_hysteresis_new}
d^c_j = \overline{d}_j \sigma_{j},\hspace{-0.75mm}
\hspace{1mm}
{\akk{\sigma_{j}^+}\icl{(t)} \hspace{-0.2mm} \in  \hspace{-0.7mm}\begin{cases}
\{1\} , \qquad \hspace{2mm} \omega_j >  \omega^1_j  \text{ or } p^c > \bar{p}^c_j, \\[1mm]
\{0\},  \qquad \hspace{2mm} \omega_j <  \omega^0_j \text{ and } p^c < \underline{p}^c_j,\\[1mm]
\{\sigma_j(t)\}, \hspace{0mm} \begin{cases}
\omega_j^0 <\hspace{-0.75mm}   \omega_j <\hspace{-0.75mm}  \omega_j^1  \hspace{-0.5mm}\text{ and } \hspace{-0.5mm} p^c \hspace{-0.75mm} < \hspace{-0.75mm}\bar{p}^c_j\hspace{-0.25mm},  \\
\omega_j \hspace{-0.75mm} < \hspace{-0.75mm} \omega_j^0 \hspace{-0.5mm}\text{ and }   \hspace{-0.5mm} \underline{p}^c_j  <\hspace{-0.75mm}  p^c <\hspace{-0.75mm}  \bar{p}^c_j,
\end{cases} \\
\{ 0, \hspace{-0.5mm} \sigma_j(t)\}, \begin{cases}
\omega_j =  \omega^0_j \text{ and }   p^c \leq \underline{p}^c_j,   \\
\omega_j \leq \omega_j^0 \text{ and }   p^c = \underline{p}^c_j,
\end{cases} \\
\{\sigma_j(t),\hspace{-0.5mm} 1\}, \begin{cases}  \omega_j = \omega^1_j \text{ and } p^c \leq \bar{p}^c_j,  \\
\omega_j \leq \omega^1_j   \text{ and } p^c = \bar{p}^c_j,
\end{cases}
\end{cases}}
\end{equation}
where $j\in N$, {$\bar{p}^c_j$ are design variables satisfying $\bar{p}^c_j > \underline{p}^c_j$,} and $\underline{p}^c_j, \overline{d}_j, \omega^0_j$ and $\omega^1_j$ are as in \eqref{sys_hysteresis_new_simple}.
Furthermore, $p^c$ follows from \eqref{sys_power_command}.
 The scheme in \eqref{sys_hysteresis_new} is depicted in Figure \ref{hysteresis_new_figure}.

\begin{figure}[t!]
\centering
\includegraphics[trim = 0mm 0mm 0mm 0mm, scale = 0.7,clip=true]{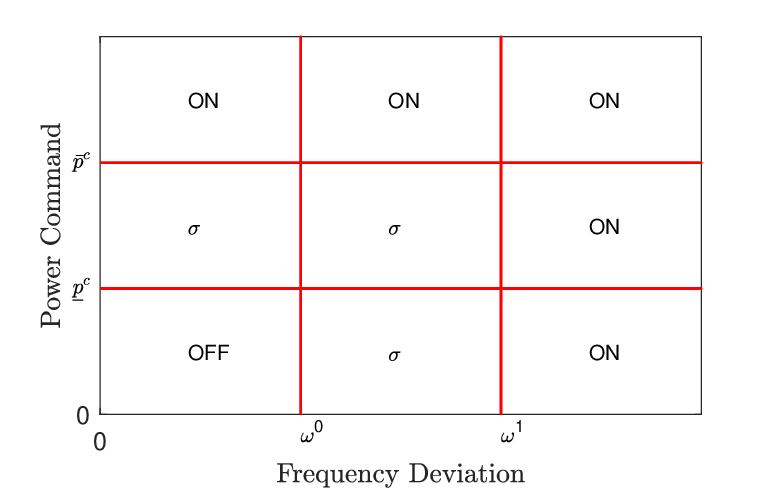}
\vspace{-0mm}
\caption{Hysteresis scheme for controllable loads described by \eqref{sys_hysteresis_new}.}
\label{hysteresis_new_figure}
\vspace{-2mm}
\end{figure}

{Compared to \eqref{sys_hysteresis_new_simple}, the scheme  in} \eqref{sys_hysteresis_new} introduces an additional threshold for power command, {such that when $p^c > \overline{p}^c_j$, then loads remain switched on.}
{As explained above, this is exploited to provide an optimality interpretation of the  resulting equilibria.}

{The H-OSLC problem \eqref{Problem_To_Min_H} is a mixed-integer optimization problem that is NP-hard \cite{karp1972reducibility}. However, the continuous relaxation of \eqref{Problem_To_Min_H} can be solved using subgradient KKT conditions (see \eqref{Problem_To_Min_RH} and Proposition \ref{prop_KKT} in \akk{the Appendix}).
Below, we {describe how}
{to appropriately design the
frequency and power command thresholds such that the  KKT conditions are satisfied by almost all loads at equilibrium, {thus leading to a power allocation that is {shown to be} very close to optimal. We also show that convergence guarantees are also provided, as in Section \ref{sec:alternative_scheme}.}
}

{To facilitate the presentation of the proposed design, let   $k \in N$ be the rank  of controllable loads when those are sorted in ascending order\footnote{
Note that in the case where there exist $i,j \in N$ such that $c^d_i/\bar{d}_i = c^d_j/\bar{d}_j$, then the order between $i$ and $j$ is {arbitrarily} assigned.}
 of $c^d_j/\bar{d}_j, j \in N$.  A parameter $x_{\underline{k}}$ is associated with the $k$-th ranked controllable load (i.e. the underlined subscript refers to the above described rank).
The design condition is presented below.}

\begin{design}\label{Design_condition}
{The values of the design variables in \eqref{sys_hysteresis_new} satisfy}
\begin{subequations}
\label{Design_cond_eqt}
\begin{gather}
\omega^0_k = c^d_k/\bar{d}_k,  k \in N,
\label{Design_cond_eqt_1} \\
\underline{p}^c_{\underline{1}} = \mathcal{D} \omega^0_{\underline{1}},
\label{Design_cond_eqt_2} \\
\underline{p}^c_{\underline{k}} = \mathcal{D} \omega^0_{\underline{k}} + \sum\limits_{j=1}^{k-1} \bar{d}_{\underline{j}}, k \in N/\{1\},
\label{Design_cond_eqt_3} \\
{\bar{p}^c_{\underline{j}} \in (\underline{p}^c_{\underline{j}}, \underline{p}^c_{\underline{j}} + \overline{\delta})}, j \in N,
\label{Design_cond_eqt_4}
\end{gather}
\end{subequations}
{where $\overline{\delta} = \min_{j \in N} \overline{d}_j$.}
\end{design}

{
\begin{remark}
Design condition \ref{Design_condition} has two important features which eliminate limit cycle behavior and also ensure that the {cost at the} resulting equilibria {is close to the optimal.}
 The choice of power command thresholds in \eqref{Design_cond_eqt_2}--\eqref{Design_cond_eqt_3} follows directly from  \eqref{sys_power_command}, \eqref{eqlb_freq}, and ensures that when $p^c \leq \underline{p}^c_j$ then $\omega^* \leq \omega^0_j$.
{This condition guarantees that loads that satisfy $p^c \leq \underline{p}^c_j$ will be switched off at equilibrium,
{which} aids in {deducing a convergence result analogous to {Theorem \ref{conv_thm_hysteresis_simple}.}}}
The optimality interpretation follows by ranking all loads based on their frequency thresholds $\omega^0_j$ and relating the latter with the cost per unit value $c^d_j / \overline{d}_j$, via \eqref{Design_cond_eqt_1}.  Then, conditions \eqref{Design_cond_eqt_2}--\eqref{Design_cond_eqt_4} ensure that when load $j$ is switched on at steady state then all loads with lower cost {per unit demand} are also switched on.
{The latter is closely linked to the KKT conditions associated with the continuous relaxation of \eqref{Problem_To_Min_H} as explained in the proof of Theorem \ref{opt_thrm_hybrid}  below.}
Condition \eqref{Design_cond_eqt_4} also ensures that $\bar{p}^c_{j} \neq \underline{p}^c_j$, {thus avoiding chattering when there is measurement noise in~$p^c$.}
\end{remark}
}

Compared to Design condition \ref{assum_pc}, Design condition \ref{Design_condition} requires knowledge of all controllable load magnitudes and also their order in terms of cost per unit $c^d_j/\bar{d}_j$, making it a centralized design.
 However, as we demonstrate in Theorem \ref{opt_thrm_hybrid} below,  Design condition \ref{Design_condition} offers {a close to optimal power allocation at steady state}. Hence, Design condition \ref{Design_condition} is preferable to Design condition \ref{assum_pc} when the required information is available. Alternatively, Design condition \ref{assum_pc} is easier to implement and requires much less information on system parameters.

{
\begin{remark}
An approach to achieve stability and optimality in power networks when on-off loads are present would be to centrally solve the mixed-integer optimization problem and then transmit the desired allocation to each load.
The scheme presented in \eqref{sys_hysteresis_new} with Design condition \ref{Design_condition} is superior to such {an} approach for two reasons. Firstly, it provides transient support to the power network, which is the main motivation for the control of on-off loads in this study.  {Secondly}, it does not require to solve \eqref{Problem_To_Min_H}, which is an NP-hard {problem} with significant computational cost when the number of loads is large.
\end{remark}
}

\subsection{Hybrid system description}

The behavior of system \eqref{sys1},~\eqref{sys2},~\eqref{sys_power_command},~\eqref{sys_hysteresis_new} can be described by the states $\zeta = (x, \sigma) \in \Lambda$.
Its continuous dynamics,
described by
\eqref{sys4_hysteresis} and \eqref{sys_power_command},
 are valid when $\zeta \in \overline{C}$ provided below.
\begin{equation}\label{f:C}
\overline{C}=\{ {\zeta} \in \Lambda:  \sigma_j \in \mathcal{J}_j(\omega_j, p^c), \;\forall j\in N\}
\end{equation}
where
\begin{equation}\label{J_set_flow}
\mathcal{J}_j(\omega_j, p^c) \hspace{-0.5mm}= \hspace{-0.5mm}\begin{cases}
\{1\}, \quad {\omega_j >  \omega^1_j,} \text{ or } p^c > \bar{p}^c_j\\[1mm]
\{0\},  \quad {
\omega_j < \omega_j^0 \text{ and }   p^c < \underline{p}^c_j, }
\\[1mm]
\{0, 1\},  \hspace{-0.5mm}\begin{cases}
\omega_j^0 \leq  \omega_j \leq \omega_j^1  \text{ and }   p^c \leq \hspace{-0.5mm}\bar{p}^c_j,  \\
\omega_j \leq \omega_j^0 \text{ and }   \underline{p}^c_j \leq p^c \leq \bar{p}^c_j.
\end{cases}
\end{cases}
\end{equation}

Alternatively, when ${\zeta}$ belongs to the set $\overline{D} ={ \Lambda \setminus \overline{C} \cup \tilde{D}}$
where
$
\tilde{D} = \{ {\zeta} \in \Lambda:
\sigma_j \in \mathcal{I}^D_j(\omega_j, p^c), \;\forall j\in N\},
$
and
\begin{equation}\label{jump_set}
\mathcal{I}^D_j(\omega_j, p^c)= \begin{cases}
\{0\}, \quad {\omega_j =  \omega^1_j \text{ or } p^c = \bar{p}^c_j},\\[1mm]
\{1 \},  \begin{cases}  \omega_j \leq \omega_j^0 \text{ and }  p^c = \underline{p}^c_j, \\
 \omega_j = \omega_j^0 \text{ and } p^c \leq \underline{p}^c_j,
\end{cases}
\end{cases}
\end{equation}
then its components follow a discrete update described by
\begin{align}\label{sys5_g}
{x^+} &= x, &
\akk{\sigma_{j}^+} = \begin{cases}
1, \quad \omega_j\geq  \omega^1_j \text{ or } p^c \geq \bar{p}^c_j,\\[1mm]
0,  \quad \omega_j\leq  \omega^0_j \text{ and } p^c \leq \underline{p}^c_j.
\end{cases}
\end{align}
 Hence, the following hybrid compact representation describes the system \eqref{sys1},~\eqref{sys2},~\eqref{sys_power_command},~\eqref{sys_hysteresis_new},
\begin{subequations}\label{sys5}
\begin{gather}
\dot{\zeta} = \overline{f}(\zeta), \zeta \in \overline{C},
\\
\zeta^+ = \overline{g}(\zeta), \zeta \in \overline{D},
\end{gather}
\end{subequations}
where $\overline{f}(\zeta): \overline{C} \rightarrow {\overline{C}}$ and $\overline{g}(\zeta): \overline{D} \rightarrow {\overline{C} \setminus \overline{D}}$  follow from \eqref{sys4_hysteresis} and \eqref{sys_power_command}, and \eqref{sys5_g} respectively.


\subsection{Analysis of equilibria and solutions}


The following proposition, proven in \akk{the Appendix}, demonstrates the existence and characterizes the equilibria of \eqref{sys5}. Note that the definition of an equilibrium to \eqref{sys5} is analogous to Definition \ref{eqbr_dfn_hybrid_2_simple} and is omitted for compactness.
\begin{proposition}\label{eqlbr_hysteresis}
Consider the system described by \eqref{sys5} and let Design condition \ref{Design_condition} hold.
Then, an equilibrium point exists and satisfies $\zeta^* \in \overline{C}$.
\end{proposition}

Proposition \ref{eqlbr_hysteresis} shows that Design condition \ref{Design_condition} suffices for the existence of equilibria to  \eqref{sys5}.
The following proposition demonstrates the existence of solutions to \eqref{sys5}, {that all maximal solutions to \eqref{sys5} are complete} and also  that no chattering occurs.

\begin{proposition}\label{Proposition_existence_no_Zeno}
For any initial condition {$\zeta(0,0) \in \Lambda$} there exists a complete solution to \eqref{sys5}. {All maximal solutions to \eqref{sys5} are complete.}
Furthermore, for any complete bounded solution to \eqref{sys5}, there exists $\tau > 0$  such that $\min_{i \geq 1} (t_{i+1, j} - t_{i, j}) \geq \tau, j \in N$.
\end{proposition}


\subsection{Stability and optimality of hybrid system}\label{sec: Stability - Optimality}

In this section, we provide our main stability and optimality results about system \eqref{sys5}, with the proofs provided in \akk{the Appendix}.

\begin{theorem}\label{conv_thm_hysteresis_opt}
Let  Design condition \ref{Design_condition} hold. Then, {for all initial conditions,}
the solutions of \eqref{sys5} {are bounded and converge} to {a subset} of its equilibria.
\end{theorem}

{Theorem \ref{conv_thm_hysteresis_opt}
and Proposition \ref{Proposition_existence_no_Zeno} demonstrate} that the inclusion of loads with dynamics described by \eqref{sys_hysteresis_new} does not compromise the stability of the {system, when  Design condition \ref{Design_condition} holds, and also does not {result to} chattering
behavior.}


{The optimality result associated with Design condition \ref{Design_condition} is  stated in Theorem \ref{opt_thrm_hybrid} below.}
 Within the {theorem statement}, we  {make} use of the {notion} of an $\epsilon$-optimal point {defined} below.

\begin{definition}\label{epsilon_close}
Given a cost function $C_f:\mathbb{R}^n \times \mathbb{Z}^m \rightarrow \mathbb{R}$ where $n,m > 0$, a vector $\bar{x} \in \mathbb{R}^n \times \mathbb{Z}^m$ is called $\epsilon$-optimal for $C_f$, {for some $\epsilon \in \mathbb{R}_{>0}$}, if it holds that
\begin{equation}
C_f(\bar{x}) \leq \min_{x \in \mathbb{R}^n \times \mathbb{Z}^m} C_f(x) + \epsilon.
\end{equation}
\end{definition}


{
\begin{theorem}\label{opt_thrm_hybrid}
Let Design condition \ref{Design_condition} hold and the control dynamics in \eqref{sys2_pm} be chosen such that $\alpha_j = c_j^{-1}, j \in N$.
 Then, the equilibrium values $(p^{M,\ast}, d^{c,\ast}, d^{u,\ast})$ are $\epsilon$-optimal for the H-OSLC problem \eqref{Problem_To_Min_H}, where
$\epsilon = \frac{1}{{2\mathcal{D}}}\max_{j \in N} (\bar{d}_j)^2$.
\end{theorem}
}

%

 Theorems  \ref{conv_thm_hysteresis_opt} and \ref{opt_thrm_hybrid} {demonstrate convergence  to {a power allocation that is close to optimal, when Design condition \ref{Design_condition} is implemented,}}
 {and provide {a non-conservative bound} on the difference between the
 {cost at equilibrium and the optimal one.}
 However, Design condition \ref{Design_condition} comes with additional information requirements compared to Design condition \ref{assum_pc}, making the latter more suitable when those parameters are difficult to obtain.

 \ak{
\begin{remark}\label{rem_epsilon_optimality}
\ic{It should be noted that the value of $\epsilon$ in Theorem \ref{opt_thrm_hybrid}, i.e. the deviation of the power allocation cost from its optimal value, is expected to be very small since in most realistic power network configurations it holds that
$\max_{j \in N} \overline{d}_j \ll \mathcal{D}$.}
\end{remark}
 }
 \begin{remark}\label{remark_exact_ell}
It should be noted that the requirement for knowledge of the  \ak{values of $\ell$ and $\mathcal{D}$ in the implementation of \eqref{sys_hysteresis_new_simple} and \eqref{sys_hysteresis_new} and Design conditions \ref{assum_pc} and \ref{Design_condition}} does not limit the applicability of the proposed schemes, since
 it can be shown that the  convergence properties presented in Theorems \ref{conv_thm_hysteresis_simple} and \ref{conv_thm_hysteresis_opt} are retained when an upper bound to \ak{the magnitude of $\ell$ and a lower bound to $\mathcal{D}$ are known.}
  However, the latter \ka{may compromise} the optimality interpretation of Theorem \ref{opt_thrm_hybrid}. Hence, there exists a  trade-off between robustness to measurement uncertainty and optimality.
\end{remark}


\ak{
\begin{remark}\label{TCNS_comparison}
The contribution of this work in comparison to \cite{kasis2017secondary_arx}, \ic{which} 
considers a hysteretic scheme 
\ic{for on-off loads as an ancillary service  to secondary frequency control},
is multilevel.
In particular, the fact that on-off loads are allowed to actively contribute at steady state when primary frequency control is considered raises issues of existence of equilibria and limit cycles.
Such issues do not occur in secondary frequency control, since loads do not contribute at equilibrium.
This study resolves these issues by providing a sufficient condition for the existence of equilibria in Theorem \ref{eqbr_hyst_existence} and a suitable hysteretic design (Design condition \ref{assum_pc}) for on-off loads which allows to deduce stability, as shown in Theorem \ref{conv_thm_hysteresis_simple}.
Moreover, using the fact that on-off loads may contribute at steady state, we proposed Design condition \ref{Design_condition}, which allows an $\epsilon$-optimal allocation among generation and on-off loads and retains the stability properties of the system, as analytically shown in Theorems  \ref{conv_thm_hysteresis_opt} and \ref{opt_thrm_hybrid}.
\end{remark}
}

\section{Simulation on the NPCC 140-bus system} \label{Simulation}

In this section we verify our {analytical} results with {numerical simulations} on the Northeast Power Coordinating
Council (NPCC) 140-bus interconnection system, using the Power System Toolbox~\cite{cheung2009power}. This model is more detailed and realistic than our analytical one, {and includes} line resistances, a DC12 exciter model{, a transient reactance generator model, and  turbine governor dynamics\footnote{The details of the simulation models can be found in the Power System Toolbox data file datanp48.}.

The test system consists of $93$ load buses serving different types of loads including constant active and reactive loads and $47$ generation buses. The overall system has a total real power of $28.55$GW. For our simulation, we added five loads on units $2, 8, 9, 16$ and $17$, each having a step increase of magnitude $3$ \ak{per unit} (base $100$MVA) at $t=1$ second.

Controllable demand was considered within the simulations on $20$ generation and $20$ load buses, with loads controlled every $10$ms.

The system was tested at three different cases. In case (i) on-off controllable loads as in \eqref{sys2_dc} were considered. The values for $\overline{\omega}_j$ were selected from a uniform distribution within the range $[0.02Hz \; 0.07Hz]$ and those of $\underline{\omega}_j$ by setting $\underline{\omega}_j = -\overline{\omega}_j$.
In case (ii) controllable loads with hysteretic dynamics  described by \eqref{sys_hysteresis} were considered. For a fair comparison, the same frequency thresholds as in case (i) were used, with $\omega^1_j = \overline{\omega}_j$ and $\omega^0_j = \omega^1_j/2$.  Finally, in case (iii), hysteretic loads following the dynamics in \eqref{sys_hysteresis_new_simple} and Design condition \ref{assum_pc} were included\footnote{Note that Design condition \ref{assum_pc} was considered in this case to provide a fair comparison between the three schemes. Later in this section, we explain how the scheme in \eqref{sys_hysteresis_new} and Design condition \ref{Design_condition} have also been implemented on the NPCC network, resulting to a stable and well behaved response.
}. For this case, the same frequency thresholds as in case (ii) where used, with power command thresholds chosen such that Design condition \ref{assum_pc} was satisfied.
For all cases $\overline{d} = 0.2$ \ak{per unit} was used.
We shall refer to cases (i), (ii), and (iii) as the 'switching', 'hysteresis' and '{adapted} hysteresis' cases respectively.

\begin{figure}[t!]
\centering
\includegraphics[trim = 1mm 0mm 0mm 0mm, scale = 0.65,clip=true]{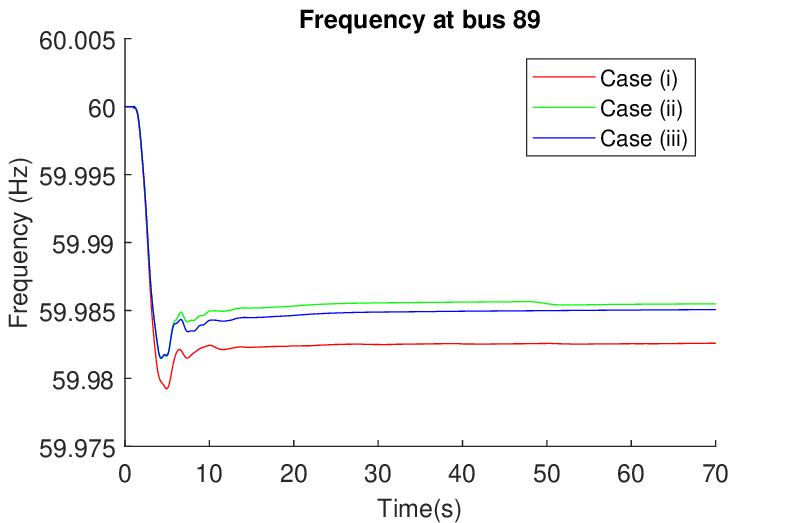}
\vspace{-2mm}
\caption{Frequency at bus 89 with controllable load dynamics as in the following three cases: i) Switching case, ii) Hysteresis case, iii) {Adapted} hysteresis case.}
\label{Frequency}
\vspace{-0mm}
\end{figure}
\begin{figure}[t!]
\centering
\includegraphics[trim = 1mm 0mm 0mm 0mm, scale = 0.65,clip=true]{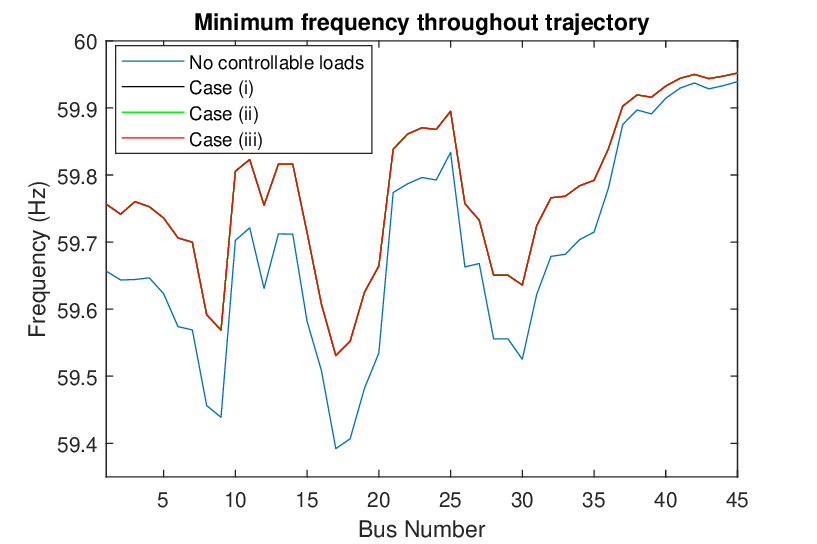}
\vspace{-2mm}
\caption{Largest frequency overshoot for buses $1-45$ for four cases: i) Switching case, ii) Hysteresis case, iii) {Adapted} hysteresis case, iv) No controllable loads case.}
\label{overshoot_freq}
\vspace{-2mm}
\end{figure}

The frequency at bus 89 for the three tested cases is  shown in Fig. \ref{Frequency}.
From this figure, we observe that {the} frequency converges to some constant value at all cases. Note that a smaller steady state frequency deviation is observed when  hysteretic loads are considered, since the {hysteresis scheme allows more loads to be switched on} at steady state {compared to \eqref{sys2_dc}}. Moreover, Fig. \ref{overshoot_freq} demonstrates that the inclusion of on-off loads decreases the maximum overshoot in frequency, by comparing the largest deviation in frequency with and without on-off controllable loads at buses $1-45$, where frequency overshoot was seen to be the largest. Note that the same {overshoot profiles are} observed in all cases (i), (ii), and (iii) since the same frequency thresholds have been used.

 Furthermore, from Fig. \ref{Fig_switch_case} it can be seen that in case (i) controllable loads switch very fast, {as demonstrated by the thick blue lines}, indicating chattering, where in case (ii) such behavior is not observed\footnote{Note that analogous behavior to case (ii) has been observed for case (iii). These results are omitted for compactness in presentation.}, since far less switches are exhibited, as shown in Fig. \ref{Fig_hysteresis_case}.
   Both figures depict the behavior at the 4 buses with hysteretic loads with the fastest consecutive switches.
{Chattering is also verified numerically in case (i),  since} it was seen that for each of the 20 controllable loads the minimum time between consecutive switches was 10ms, which is the smallest time increment in our discrete numerical simulation.
Therefore, the numerical results support the analysis of this paper, verifying that hysteresis eliminates chattering of controllable loads.

To demonstrate the possibility of limit cycles when case (ii) is considered, we altered the {frequency thresholds} of the on-off load at bus $21$, making $\omega^1_{21}$ coincide with the equilibrium frequency and then repeated the simulations for cases (ii) and (iii). As demonstrated on Fig. \ref{Limit_cycle_case}, the load at bus $21$ exhibits limit cycle behavior at steady state, whereas when the {adapted} hysteresis scheme was considered, no such behavior is observed.

\begin{figure}[t!]
\centering
\includegraphics[trim = 1mm 0mm 0mm 0mm,  scale = 0.6, clip=true]{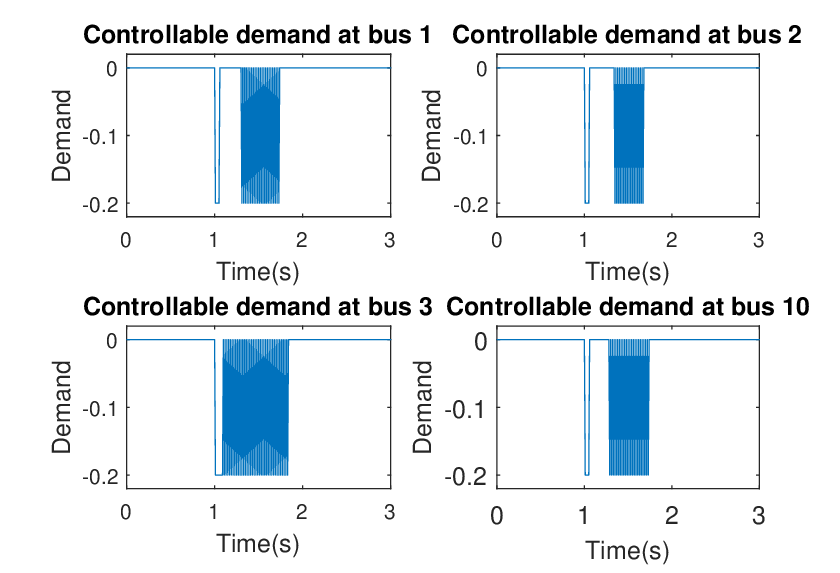}
\vspace{-2.7mm}
\caption{Controllable demand at 4 buses with on-off loads described by \eqref{sys2_dc}.}
\label{Fig_switch_case}
\vspace{-0mm}
\end{figure}

\begin{figure}[t!]
\centering
\includegraphics[trim = 1mm 0mm 0mm 0mm, scale = 0.6,clip=true]{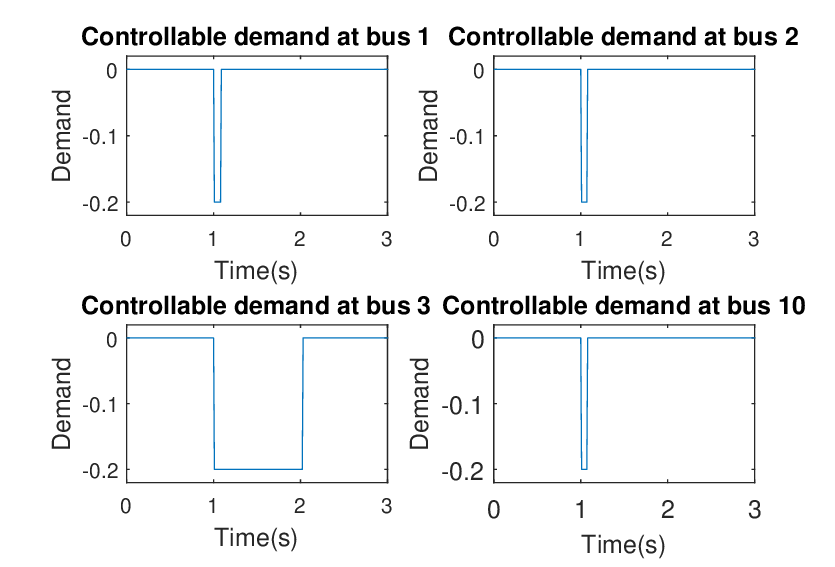}
\vspace{-0mm}
\caption{Controllable demand at 4 buses with Hysteretic on-off loads.}
\label{Fig_hysteresis_case}
\vspace{-0mm}
\end{figure}

\begin{figure}[t!]
\centering
\includegraphics[trim = 1mm 0mm 0mm 0mm,  scale = 0.65,clip=true]{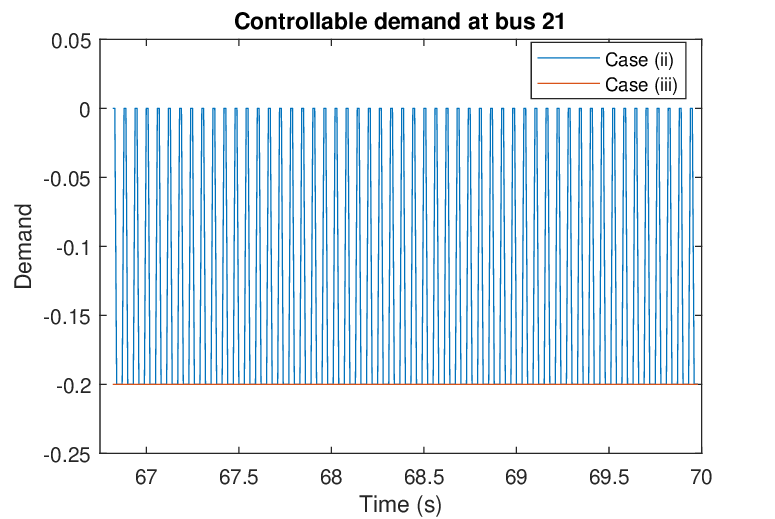}
\vspace{-0mm}
\caption{Controllable demand at bus $21$ for cases (ii) and (iii).}
\label{Limit_cycle_case}
\vspace{-2mm}
\end{figure}

To verify the optimality results of Theorem \ref{opt_thrm_hybrid}, we repeated the simulation with 47 loads on generation buses with magnitudes randomly selected from a uniform distribution of range $[0.025\; 0.075]$ \ak{per unit} and 20 loads on load buses $1-20$ of magnitude 0.2 \ak{per unit}. We aimed for a larger number of on-off devices to allow a large number of {possible solutions to \eqref{Problem_To_Min_H}} and show that the cost of the obtained equilibrium when Design condition \ref{Design_condition} is applied is $\epsilon$-close to the globally minimum of \eqref{Problem_To_Min_H}, as follows from Theorem \ref{opt_thrm_hybrid}. The costs were selected from a uniform distribution in the range $[10^{-4} \; 0.002]$, resulting to frequency thresholds $\omega^0_j$ following \eqref{Design_cond_eqt_1} in the range $[0.013Hz \; 0.08Hz]$.
 The simulation results verified the convergence of frequency to an equilibrium value and that no limit cycle behavior occurred, similarly to case (iii). Furthermore, the obtained equilibrium was seen to be identical with the optimal one, calculated using a heuristic reproduction and mutation genetic algorithm (see \cite[Ch. 3]{goldberg1988genetic}) implemented 1000 times with random initial conditions and always converging to the same solution, numerically verifying the optimality analysis of Theorem \ref{opt_thrm_hybrid} when Design condition \ref{Design_condition} is considered.

\section{Conclusion}\label{Conclusion}
We have considered the problem of primary frequency control where controllable on-off loads provide ancillary services to the power network.  We first considered loads that switch on when some frequency threshold is reached and off otherwise and provided relevant stability guarantees for the power network. Furthermore, {we} discussed that such schemes might exhibit {chattering, which {limits} their practicality. To cope with this, we considered on-off loads with hysteretic dynamics and {showed} that  chattering is no longer exhibited.} {We also} provided {design} conditions that guarantee the existence of equilibria when such loads are considered.
 However, numerical simulations demonstrate that hysteretic loads may exhibit limit cycle behavior.
  As a remedy to this problem, we proposed an adapted {hysteretic} control scheme and appropriate design conditions
  {that ensure that the network is stable, while also avoiding chattering.}
Furthermore, we considered {a mixed-integer optimization problem for power allocation. We} {proposed a suitable control design such that {the} stability guarantees are retained and  the {cost at the equilibria of the system is within}} $\epsilon$ to the global {minimum, providing also} a non-conservative bound for~$\epsilon$.
Our {analytical} results have been verified with numerical simulations on the NPCC 140-bus system where it was shown that the presence of on-off loads reduces the frequency overshoot and that {hysteretic} schemes {avoid} chattering. Furthermore, simulation results demonstrate that our proposed hysteretic scheme {avoids} 
limit cycle behavior and {leads} to an optimal power allocation at steady state.

\vspace{-0.0mm}
\section*{Appendix}
\vspace{-0.0mm}

This appendix {includes} the proofs of {the results presented in the main text}.

\emph{Proof of Lemma~\ref{Uniqueness_Existence}:}
{The lemma can be proven using Proposition 3 in \cite{cortes2008discontinuous}. This states that solutions exist if $Q$
 \ak{is locally essentially bounded and measurable.
 The local boundedness of $Q$ follows since the step size at discontinuities is always bounded from \eqref{sys2_dc} (i.e. by $\max_{j \in N} \overline{d}_j$)
   and the Lipschitz property of the \ic{vector \ka{field} in} \eqref{sys1}, \eqref{sys2}.
Moreover, the fact that $Q$ is measurable follows trivially.}
    \hfill $\blacksquare$

Within the proof of Theorem \ref{convthm} we will make use of the following equilibrium equations for system \eqref{sys1}--\eqref{sys2_dc}, which follow from {Definition \ref{eqbrdef}}.
\begin{subequations} \label{eqbr}
\begin{gather}
0 = \omega^*_i - \omega^*_j, \; (i,j) \in E, \label{eqbr1} \\
0 \in - p_j^L \! + \! p_j^{M,*} \! -F[d^c_j] - \!\! \sum_{\ak{k \in N^o_j}} p^*_{jk} + \!\! \sum_{\ak{i \in N^i_j}} p^*_{ij}, \; j\in N, \label{eqbr2} \\
p^*_{ij}=B_{ij} \eta^*_{ij}, \; (i,j) \in E,\ \label{eqbr4} \\
p^{M,*}_j = -\alpha_j\omega^*_j, \;  d^{u,*}_j = A_j \omega^*_j,  j \in N. \label{eqbr5}
\end{gather}
\end{subequations}
{The notion of a Lyapunov stable equilibrium point will also be used for the discontinuous system considered and is defined below.}
{\begin{definition}\label{Lyap_stable}
An equilibrium point $x^\ast$ of \eqref{sys_Filippov_representation} is Lyapunov stable if for all $\epsilon>0$ there exists a $\delta>0$ s.t. any Filippov solution $x(t)$ of  \eqref{sys1}-\eqref{sys2_dc} with initial condition $x(0)=x_0$, $\|x_0-x^\ast\|<\delta$, satisfies $\|x(t)-x^\ast\|<\epsilon$ for all $t\geq0$.
\end{definition}}

\emph{Proof of Theorem~\ref{convthm}:}
{To prove Theorem~\ref{convthm}, we {make} use of \cite[Theorem 3]{bacciotti1999stability} to establish convergence to the set of equilibria of \eqref{sys_Filippov_representation}. We {then} show that each equilibrium to \eqref{sys_Filippov_representation} is Lyapunov stable and deduce convergence to an equilibrium point using similar arguments as in  \cite[Prop. 4.7, Thm. 4.20]{haddad2011nonlinear}.}

We will use the dynamics in~\eqref{sys1} and \eqref{sys2} to define a Lyapunov function for system~\eqref{sys1}--\eqref{sys2_dc}. {Note that the set valued map in \eqref{sys_Filippov_representation} takes compact, convex values, in accordance to the class of systems considered in \cite{bacciotti1999stability}.}

Firstly, we consider {some equilibrium point $x^\ast = (\eta^*, \omega^*, p^{M,*})$ and the function} $V_F (\omega) = \frac{1}{2}\sum_{j \in N} M_j (\omega_j- \omega^*_j)^2$. {We {then} consider the time-derivative of $V_F (\omega)$ along the solutions of~\eqref{sys1}--\eqref{sys2_dc}. For a given value of the state $x=(\eta, \omega, p^{M})$ this is a set valued map given by} {${\dot{V}_F\ak{(x)} := \{\frac{\partial V_F}{\partial x}\dot x: \dot x \in Q(x)\} = }\{\sum_{j \in N} (\omega_j - \omega^*_j) (-p^L_j + p_j^M - u_j - d^u_j - \sum_{\ak{k \in N^0_j}} p_{jk} + \sum_{\ak{i \in N^i_j}} p_{ij}): u_j \in F[d^c_j(\omega_j)]\}$}, by substituting~\eqref{sys1b} for} $\dot{\omega}_j$  and using the differential inclusion for $d^c_j$ for $j \in N$. Subtracting the product of $(\omega_j - \omega^*_j)$ with each term in~\eqref{eqbr2}, we get
\begin{align}
&\dot{V}_F\ak{(x)} \hspace{-0.75mm} =\hspace{-0.75mm}  \{\hspace{-0.5mm} \sum_{j \in N}\hspace{-0.5mm} (\omega_j \hspace{-0.75mm}- \hspace{-0.75mm}\omega^*_j) (p^M_j \hspace{-0.5mm}-\hspace{-0.5mm} p^{M,*}_j \hspace{-0.75mm}-(u_j \hspace{-0.75mm}-\hspace{-0.75mm} u^*_j)\hspace{-0.5mm} -\hspace{-0.5mm} (d^u_j \hspace{-0.75mm}-\hspace{-0.75mm} d^{u,*}_j))
\nonumber\\ & + \hspace{-2.5mm} \sum_{(i,j) \in E} \hspace{-2.0mm}(p_{ij} \hspace{-0.5mm}-\hspace{-0.5mm} p^*_{ij}) (\omega_j \hspace{-0.5mm}-\hspace{-0.5mm} \omega_i)
\hspace{-0.5mm}:\hspace{-0.5mm} u_j \hspace{-0.5mm} \in \hspace{-0.5mm} F[d^c_j(\omega_j)],\hspace{-0.5mm} u^*_j\hspace{-0.5mm} \in \hspace{-0.5mm}F[d^c_j(\omega^*_j)]\}, \label{VFdiff}
\end{align}
using in the first term the equilibrium {condition~\eqref{eqbr1}.}

Additionally, consider $V_P(\eta) = \sum_{(i,j) \in E} B_{ij} (\eta_{ij} - \eta_{ij}^*)^2$.
 Using~\eqref{sys1a} and~\eqref{sys1d}, the time-derivative equals
\begin{align}
\hspace{-0.5mm}\dot{V}_P\ka{(x)} &\hspace{-0.5mm}= \hspace{-2.5mm}\sum_{(i,j) \in E}\hspace{-2.5mm} B_{ij} (\eta_{ij} \hspace{-0.5mm}-\hspace{-0.5mm} \eta^*_{ij}) (\omega_i \hspace{-0.5mm}-\hspace{-0.5mm} \omega_j)
\hspace{-0.5mm} \nonumber \\
&= \hspace{-2.5mm}\sum_{(i,j) \in E}\hspace{-2.0mm} (p_{ij} \hspace{-0.5mm}-\hspace{-0.5mm} p^*_{ij}) (\omega_i \hspace{-0.5mm}-\hspace{-0.5mm} \omega_j). \label{VPdiff}
\end{align}
Finally, consider the function $V_M(p^M) = \frac{1}{2}\sum_{j \in N} \tau_j(p^M_j - p^{M,*}_j)^2$. Using \eqref{sys2_pm}, its time derivative is given by
\begin{equation}\label{VMdiff}
\dot{V}_M\ka{(x)}  \hspace{-0.5mm}= \hspace{-1mm}\sum_{j \in N} (p^M_j \hspace{-0.5mm}-\hspace{-0.5mm} p^{M,*}_j)(-(p^M_j \hspace{-0.5mm}-\hspace{-0.5mm} p^{M,*}_j) \hspace{-0.5mm}-\hspace{-0.5mm} (\omega_j \hspace{-0.5mm}-\hspace{-0.5mm} \omega^*_j)).
\end{equation}

Based on the above, we define the function
\begin{align}\label{V_Lyapuinov}
&V(\eta, \omega, p^M) = V_F(\omega) + V_P(\eta) + V_M(p^M),
\end{align}
which is continuously differentiable and {has a strict minimum at $(\eta^*, \omega^*, p^{M,*})$} and hence is a suitable Lyapunov candidate. {Furthermore, it trivially follows that $V$ is regular, following the definition provided in  \cite[p.363-364]{bacciotti1999stability}.}
By \eqref{sys2_du} and~\eqref{VFdiff}--\eqref{VMdiff}, it follows that
\begin{align*}
& \dot{V}\ak{(x)} = \{\sum_{j \in N} [ - A_j(\omega_j - \omega^*_j)^2 - (p^M_j - p^{M,*}_j)^2 \nonumber \\
& -(\omega_j - \omega^*_j) (u_j - u^*_j)] : u_j \in F[d^c_j(\omega_j)], u^*_j \in F[d^c_j(\omega^*_j)]\}.
 \end{align*}
  Using~\eqref{sys2_dc_set_valued}, it therefore holds that,
\begin{align}
{\max_{y \in\dot{V}\ak{(x)}} y  }& \ \le \hspace{-0.5mm}\sum_{j \in N} \hspace{-0.5mm}[-A_j(\omega_j - \omega^*_j)^2 - (p^M_j - p^{M,*}_j)^2]
\le 0. \label{vdotineq}
\end{align}

It is clear that $V\ak{(x)}$ has a global minimum at $x^* = (\eta^*, \omega^*, p^{M,*})$.
Furthermore, from \eqref{vdotineq}, there exists {\ic{a compact set  $\Xi = \{x \colon V(x) \ic{- V(x^\ast)} \le \epsilon
 \},$ for some $\epsilon > 0$,}}
such that solutions initiated in $\Xi$ remain in $\Xi$ for all future times.
Note that the value of $\epsilon$ in the definition of $\Xi$ can be selected to be arbitrarily small and hence $x^\ast$ is Lyapunov stable, following Definition \ref{Lyap_stable}. Moreover, note that $x^\ast$ is an arbitrarily selected equilibrium point, which allows to extend the above argument to all equilibria of \eqref{sys_Filippov_representation}.
}

Therefore, Theorem 3 in \cite{bacciotti1999stability} can be used  on function $V\ak{(\eta, \omega, p^M)}$ within the {compact set} $\Xi$ along solutions of~\eqref{sys1}--\eqref{sys2_dc}.
{Let $Z = \{\ic{x} : 0 \in \dot{V}\ak{(x)}\}$ and $\Psi$ be the largest weakly\footnote{The notion of a weakly invariant set used follows from \cite[Dfn. 4]{bacciotti1999stability}. } invariant set within $\Xi \cap Z$. Then, Theorem 3 in \cite{bacciotti1999stability}  guarantees that all solutions of ~\eqref{sys1}--\eqref{sys2_dc}}
 that start within $\Xi$ converge to $\Psi$.
 \ak{Since $0 \in \dot{V}\ak{(x)}$,
it follows that within $\Psi$, $\omega = \omega^*$ and $p^M = p^{M,*}$ from \eqref{vdotineq}.
Moreover,  \eqref{sys1a} suggests that when $\omega$ is constant then  $\eta$ also takes some constant value $\bar{\eta}$.}
Applying the above to equation~\eqref{sys1},  leads to the equilibrium conditions~\eqref{eqbr1}--\eqref{eqbr5}. Therefore, we conclude by Theorem 3 in  \cite{bacciotti1999stability} that all Filippov solutions of~\eqref{sys1}--\eqref{sys2_dc} with initial conditions $(\eta(0), \omega(0), p^M(0)) \in \Xi$ converge to the set of equilibria {within $\Xi$} defined in {Definition~\ref{eqbrdef}.}

{Since $\Xi$ is a bounded set for given $\epsilon$ and solutions of  \eqref{sys_Filippov_representation} initiated within $\Xi$ converge to the set of its equilibria, it follows that trajectories of \eqref{sys_Filippov_representation} are always bounded and hence each trajectory $x(t)$ has an $\omega$-limit point (see  \cite[p.129]{filippov1988differential})  that is an equilibrium, i.e. there exists a subsequence $x(t_n)$ that converges to an equilibrium point as $t \rightarrow \infty$. Since all equilibria are also Lyapunov stable, it follows that each trajectory initiated within $\Xi$ converges to an equilibrium within $\Xi$. Hence, noting that $\Xi$ can be {chosen to be} arbitrarily large, we deduce global convergence of solutions to \eqref{sys_Filippov_representation} to an equilibrium point of \eqref{sys_Filippov_representation}, which completes the proof.
}
 \hfill $\blacksquare$

The following lemma characterizes the equilibria of \eqref{sys4} and will be used in the proof of Theorem \ref{eqbr_hyst_existence}.
\begin{lemma}\label{eqlbr_hysteresis_1}
Consider the system described by \eqref{sys4}.
Then any equilibrium point $\zeta^* = (x^*, \sigma^*) \in C$.
\end{lemma}
\emph{Proof of Lemma~\ref{eqlbr_hysteresis_1}:}
Recall from Definition \ref{eqbr_dfn_hybrid} that any equilibrium $\zeta^*$ should satisfy $f(\zeta^*) = 0, \zeta^* \in C$ or $\zeta^+ = \zeta^*, \zeta^* \in D$. Now note that the latter case can be excluded since $g(\zeta):D \rightarrow {C \setminus D}$. Therefore, all $\zeta^* \in C$.
\hfill $\blacksquare$

\emph{Proof of Theorem \ref{eqbr_hyst_existence}:}  Within the proof we use
 $h_j = \omega^1_j - \omega^0_j$ and
\begin{equation}\label{omega_star}
\tilde{\omega}(p^L, \sigma) = (-\ell - \overline{d}^T \sigma) / \mathcal{D},
\end{equation}
reminding that $\ell = \vect{1}^T_{|N|} p^L$.
 As follows from \eqref{sys1b} and Lemma \ref{eqlbr_hysteresis_1}, the existence of an equilibrium to \eqref{sys4} is equivalent {to} the existence of a pair $(\tilde{\omega}, \tilde{\sigma})$ within some equilibrium point $\zeta^* \in C$ as follows from Definition \ref{eqbr_dfn_hybrid}, such that \eqref{omega_star} is satisfied.
 Furthermore,  for scalar $\tilde{\omega}_i$ and vector $\tilde{\sigma}$, we define the set {$\Pi$} as
$\Pi(\tilde{\omega}_i, \tilde{\sigma}) = \{ j: \tilde{\omega}_i > \omega^1_j \text{ and } {\tilde{\sigma_j} = 0, \tilde{\omega}_i < \omega^0_j \text{ and } \tilde{\sigma_j} = 1\}}$, i.e. it contains all buses that violate \eqref{sys_hysteresis} when $\omega = \tilde{\omega}_i \vect{1}_{|N|}$ and $\sigma = \tilde{\sigma}$. It should be clear that for any feasible equilibrium with frequency $\omega^*$ and switching state $\sigma^*$, then ${\Pi}(\tilde{\omega}_i = \omega^*,  \sigma^*) = \emptyset$.


To show that the condition {suffices for the existence of equilibria}, we prove that when $h_i \geq \overline{d}_i/\mathcal{D}$, then there exists some $\zeta^*$ that satisfies Definition \ref{eqbr_dfn_hybrid} for any $p^L$. An  equilibrium pair $(\omega^*, \sigma^*)$  may be obtained by the following algorithm.
Consider any $p^L$, a vector $\sigma^0 =  \vect{0}_{|N|}$ and the corresponding value of  $\tilde{\omega}_0 = \tilde{\omega}(p^L,\sigma^0)$, as follows from
\eqref{omega_star}.
Then consider the set ${\Pi_0 = \Pi}(\tilde{\omega}_0, \sigma^0)$
and note that it contains all the buses with $\sigma_j = 0$ that should satisfy $\sigma_j = 1$ when $\omega_j = \tilde{\omega}_0$, as follows from \eqref{sys_hysteresis}.
Then, choose the bus $j$ that satisfies
 $\omega^0_j = \min_{k \in \Pi_0} \omega^0_k$ and
 define $\sigma^1 = \{ \sigma : \sigma_i = \sigma^0_i, i \in N/\{j\}, \sigma_j = 1\}$,
 Then,  $\tilde{\omega}_1 = \tilde{\omega}_0 - \overline{d}_j/\mathcal{D}$, noting that the condition $h_i \geq \overline{d}_i/\mathcal{D}, i \in N$ guarantees that $\tilde{\omega}_1 > \omega^0_j$. Then, define the set ${\Pi_1 = \Pi(\tilde{\omega}_1, \sigma^1})$ and repeat. This algorithm creates a decreasing series of $\tilde{\omega}_i$ and a series of {$\sigma^i$} that converge to some values $(\omega^*, \sigma^*)$ after at most $|\Pi_0|$ iterations. This follows, by noting that it always holds that when $\Pi_i \neq \emptyset$, then $|\Pi_{i+1}| \leq |\Pi_i| - 1$, since no bus with $\sigma_j = 1$ belongs to any set $\Pi_i$, since at any iteration $\tilde{\omega}_i > \omega^0_j$ for any $j$ where $\sigma_j =1$.
Hence, the algorithm converges after at most $|\Pi_0|$ iterations to some pair ($\omega^*, \sigma^*$) that satisfies both  \eqref{sys_hysteresis} and \eqref{omega_star}.
Therefore,  when $h_i \geq \overline{d}_i/\mathcal{D}$ there exists an equilibrium that satisfies Definition \ref{eqbr_dfn_hybrid}.
 \hfill $\blacksquare$

\emph{Proofs of Propositions~\ref{dwell_time_lemma},~\ref{Proposition_existence_no_Zeno_simple} and~\ref{Proposition_existence_no_Zeno}:}
All proofs follow in analogy to the proofs of Lemma 4 and Proposition 1 in \cite{kasis2017secondary_arx}. {Note that the fact that all maximal solutions are complete follows from the global Lipschitz properties of the continuous variable $x$ in \eqref{sys4}, \eqref{sys5_simple} and \eqref{sys5}.}
\hfill $\blacksquare$

\emph{Proof of Proposition~\ref{eqlbr_hysteresis_simple}:}
Recall from Definition \ref{eqbr_dfn_hybrid_2_simple} that any equilibrium {$\zeta^*$} should satisfy {$\tilde{f}(\zeta^*) = 0, \zeta^* \in F$ or $\zeta^+ = \zeta^*, \zeta^* \in G$}. Now note that the latter case can be excluded since ${\tilde{g}(\zeta)}:G \rightarrow {F \setminus G}$. Therefore, all ${\zeta^*} \in F$. To show that an equilibrium of \eqref{sys5_simple} exists, it suffices that \eqref{sys_hysteresis_new_simple},~\eqref{sys_power_command} and \eqref{eqlb_freq} are simultaneously satisfied for some $\omega$ and $\sigma$. Now define the set of buses $N_1 = \{j : p^c \leq \underline{p}^c_j\}$.
 Then, there exists an equilibrium with $\sigma_j = 0, j \in N_1$ and $\sigma_j = 1, j \in N \setminus N_1$ that satisfies  \eqref{sys_hysteresis_new_simple}, \eqref{sys_power_command} and \eqref{eqlb_freq}, when  Design condition \ref{assum_pc}  holds.
\hfill $\blacksquare$


\emph{Proof of Proposition~\ref{eqlbr_hysteresis}:}
The proof follows analogously to the proof of Proposition~\ref{eqlbr_hysteresis_simple}, noting that the constructed equilibrium in the last argument is also in agreement with Design condition \ref{Design_condition}.
\hfill $\blacksquare$


{Note the for convenience we prove first Theorem~\ref{conv_thm_hysteresis_opt}, before proving Theorem~\ref{conv_thm_hysteresis_simple}, as the latter follows easily from the proof of the former.} {Furthermore note that both systems \eqref{sys5_simple} and \eqref{sys5} considered in Theorems \ref{conv_thm_hysteresis_simple} and \ref{conv_thm_hysteresis_opt} respectively are well posed, satisfying the conditions in \cite[Theorem 6.8]{goebel2012hybrid}.}

\emph{Proof of Theorem~\ref{conv_thm_hysteresis_opt}:}
%
{To prove Theorem \ref{conv_thm_hysteresis_opt} we
 first}
define  the sets of buses {$N_1 = \{j : p^c \leq \underline{p}^c_j\}$} and $N_2 = N \setminus N_1$.
{We then split the proof in two parts. In part (a), we}
 {show that {for each initial condition} there exists some finite time $T$ such that for each $j \in N_2$ it either holds that (i) $\sigma_j(t) = \sigma^*_j, t \geq T$ or (ii) $\omega^* = \omega^1_j$ {for all solutions to \eqref{sys5}.
 Then, {in part (b)}  we show that {when either of these two cases holds, a Lyapuonv argument (\cite[Corollary 8.7 (b)]{goebel2012hybrid}) can be used to deduce convergence to the  set of equilibria of \eqref{sys4}. }}}

{\textbf{Part (a):}} When Design condition \ref{Design_condition} holds {the equilibrium frequency $\omega^*$  satisfies} $\omega^* \leq \omega^0_j, j \in N_1$ .
This follows
from the monotonicity in the map from $\ell$ to $p^c$ and the fact that when $p^c = \underline{p}^c_k$ as follows from \eqref{Design_cond_eqt_2}--\eqref{Design_cond_eqt_3}, then $\omega^* = \omega^0_k$, as follows from \eqref{sys_power_command},~\eqref{eqlb_freq},~\eqref{Design_cond_eqt}.
Furthermore, {for each initial condition}, when $j \in N_2$, it holds that either $\sigma_j$ converges to some $\sigma_j^*$ in some finite time
 $T$ following the fact that $\sigma_j$ is not allowed to switch from $1$ to $0$ from \eqref{sys_hysteresis_new} or {that}
 }}
  $\omega^* = \omega^1_j$.
 Note that the two above cases
 are not mutually exclusive.

{\textbf{Part (b):}} {In this part, we show that when \eqref{sys5} satisfies
 either  (i) $\sigma_j(t) = \sigma^*_j, t \geq T$ or (ii) $\omega^* = \omega^1_j$
for $j \in N_2$, then a Lyapunov argument can be used to show that for all initial conditions at time $T$ solutions to \eqref{sys5} convergence to a subset of its equilibria.}
{First, consider the} continuous function $V$, described by \eqref{V_Lyapuinov}.
 Using similar arguments as in the proof of Theorem \ref{convthm} and defining $T_c =  \{{t  \geq T}: (t,{\ell}) \in K, {\zeta}(t,{\ell}) \in \overline{C} \}, T_d =  \{{t  \geq T}: (t,{\ell}) \in K, {\zeta}(t,{\ell}) \in \overline{D} \}$, {where $K$ is a hybrid time domain for \eqref{sys5} and $\overline{C}$ and $\overline{D}$ are defined with the aid of \eqref{f:C} and \eqref{jump_set} respectively,} it follows that
\begin{subequations}
\begin{gather}
\dot{V}\ak{(x)} \hspace{-0mm} = \hspace{-0mm}-\hspace{-0mm} \sum_{j \in N}\hspace{-0.0mm} [\hspace{-0mm}A_j(\omega_j \hspace{-0mm}-\hspace{-0mm} \omega^*_j)^2 \hspace{-0mm}+\hspace{-0mm} (p^M_j \hspace{-0mm}-\hspace{-0mm} p^{M,*}_j)^2 \nonumber \\
\hspace{-0mm}+\hspace{-0mm}(\omega_j \hspace{-0mm}-\hspace{-0mm} \omega^*_j)(d^c_j \hspace{-0mm}-\hspace{-0mm} d^{c,*}_j) ]
\nonumber \\
\leq  \sum_{j \in N} [-A_j(\omega_j - \omega^*_j)^2 - (p^M_j - p^{M,*}_j)^2] ,  t \in T_c \label{Lyap_continuous_derivative}
\\
V(\zeta^+) - V(\zeta) = 0, t \in T_d,
\end{gather}
\end{subequations}
along any solution of \eqref{sys5}, {where $\zeta^+ = (x^+, \sigma^+)$.
{Note that $(\omega_j - \omega^*_j)(d^c_j - d^{c,*}_j) \geq 0$ in \eqref{Lyap_continuous_derivative} follows  since $\omega^*~\leq~\omega^0_j, j \in N_1$ {(shown in part(a))}, and the fact that {in part (b) we consider that} for $t \geq T$ it either holds that $\sigma_j = \sigma^*_j$ or $\omega^* = \omega^1_j$ for {$j \in N_2$}. 
}
 Furthermore, note that when ${x \subset \zeta \in \overline{D}}$, the value of $V\ak{(x)}$ remains constant as it only depends on $x$ that is constant from \eqref{sys5_g}.
Note that $V\ak{(x)}$ has a strict minimum at $(\eta^*, \omega^*, p^{M,*})$.
Moreover, $V\ak{(x)}=0$ yields $(\eta, \omega, p^M) =(\eta^\ast, \omega^\ast, p^{M,\ast})$, and thus $\sigma=\sigma^\ast$. Hence, the function $V$ serves as a Lyapunov function for the hybrid system \eqref{sys5}. Then, there exists a  set
$
S=\{(x, \sigma): x\in \Xi {\text{ and }} \sigma\in {\mathcal{J}}(\omega, p^c)\}$ for some neighborhood $\Xi$ of $x^\ast$,
where
\ka{$\Xi$ is a compact set  satisfying} \ka{$\Xi = \{x \colon V(x) \ic{- V(x^\ast)} \le \epsilon
 \},$ for some $\epsilon > 0$,}
{such that solutions to \eqref{sys5} that lie in $S$ at $t = T$, stay within $S$ for all future times.}
{Moreover, note that the set $\Xi$ is compact, and hence solutions within $S$ are bounded.}
Furthermore, as shown in Proposition \ref{Proposition_existence_no_Zeno}, {all maximal solutions to \eqref{sys5} are complete,}
 and {for bounded solutions to \eqref{sys5},} the time interval between any two consecutive switches {of individual loads} is lower bounded by a positive number.
Therefore, by \cite[Corollary 8.7 (b)]{goebel2012hybrid},  there exists $r>0$ such that
{all complete and bounded}
solutions to~\eqref{sys5} {with initial conditions at time $T$ in $S$} converge to the largest weakly invariant\footnote{The definition of a weakly invariant set to a hybrid system is provided in \cite[Dfn. 6.19]{goebel2012hybrid}.}
 subset
 of the set
$\{ \zeta : V(x) = r\} \cap \{{\zeta:\zeta \in \overline{C},} \dot{V}\ak{(x)} = 0\} \cap S$, {which corresponds to a set of equilibria within $S$}.
The characterization of this invariant set and the fact that $\Xi$ can be arbitrarily large follows in a similar way as in the proof of Theorem \ref{convthm}, noting that the equilibria of \eqref{sys5} are as described by Proposition \ref{eqlbr_hysteresis}.

{Noting that in part (a) we showed that for each initial condition there exists a time $T$ such that for each $j \in N_2$ either  (i) $\sigma_j(t) = \sigma^*_j, t \geq T$ or (ii) $\omega^* = \omega^1_j$ {holds,} allows to deduce Theorem \ref{conv_thm_hysteresis_opt} and completes the proof.}
\hfill $\blacksquare$

\emph{Proof of Theorem~\ref{conv_thm_hysteresis_simple}:}
{
In analogy to the proof of Theorem \ref{conv_thm_hysteresis_opt},  note that, when Design condition \ref{assum_pc} holds,    $\omega^* \leq \omega^0_j, j \in N_1$. The latter follows directly from the  equations for power command and equilibrium frequency, described in \eqref{sys_power_command},~\eqref{eqlb_freq}. Alternatively, when $j \in N_2$, the same arguments as when Design condition \ref{Design_condition} is considered hold, to deduce that that either $\sigma_j$ converges to some $\sigma_j^*, j \in N_2$ in some finite time
 $T$ following the fact that $\sigma_j$ is not allowed to switch from $1$ to $0$ from \eqref{sys_hysteresis_new_simple}, or that $\omega^* = \omega^1_j$.
 The rest of the proof follows analogously to the proof of Theorem~\ref{conv_thm_hysteresis_opt}.
 This is since, for given $\ell$, the {loads} that satisfy $p^c \geq \overline{p}^c_j$ {have constant demand} and hence {those} do not affect the dynamic behavior of
\eqref{sys_hysteresis_new},  which reduces to that of \eqref{sys_hysteresis_new_simple}.
}
\hfill $\blacksquare$


{Within the} the proof of Theorem \ref{opt_thrm_hybrid}, we  {consider} a relaxed version of the H-OSLC problem \eqref{Problem_To_Min_H} by allowing continuous values for controllable loads.
Furthermore, we  {relax} the  discrete cost functions $C_{h,j}$ to $\hat{C}_{h,j}$ as follows:
\begin{equation}
\hat{C}_{h,j}(d^c_j) = \begin{cases}
\gamma_j d^c_j, \; 0 \leq d^c_j \leq \bar{d}_j, \\
\infty, \text{ otherwise},
\end{cases}  j \in N,
\label{CD_continuous}
\end{equation}
where $\gamma_j = c^d_j / \bar{d}_j$.  Hence, we define the following optimization problem, called the relaxed hybrid optimal supply and load control problem (RH - OSLC)
\begin{equation}
\begin{aligned}
&\hspace{2em}\underline{\text{RH - OSLC:}} \\
&\min_{p^M,d^c,d^u} \sum\limits_{j\in N} \Big( \frac{c_j}{2}(p_{j}^M)^2 +   \hat{C}_{h,j} (d^c_{j}) + \frac{1}{2A_j}(d^{u}_j)^2  \Big) \\
&\text{subject to } \sum\limits_{j\in N} (p_j^M - d^u_j - p_j^L) =  \sum\limits_{j\in  N} d^c_j.
 \label{Problem_To_Min_RH}
\end{aligned}
\end{equation}

The RH-OSLC problem is convex since each component of the cost function is convex. To solve the $RH-OSLC$ problem we shall make use of subgradient techniques \cite[Section 23]{rockafellar2015convex} {and} the KKT conditions, as follows from Proposition \ref{prop_KKT} below,
where $\partial \hat{C}_{h,j}(\bar{d}^c_j)$ denotes the subdifferential of $\hat{C}_{h,j}$ at $\bar{d}^c_j$ (see e.g. \cite{rockafellar2015convex}).

\begin{proposition}\label{prop_KKT}
A point $(\bar{p}^M, \bar{d}^c, \bar{d}^u)$ is a global minimum of \eqref{Problem_To_Min_RH} if and only if there exists $\lambda \in \mathbb{R}$ such that
\begin{subequations}\label{kkt_rh}
\begin{equation}\label{kkt1_rh}
\sum\limits_{j\in  N} (\bar{d}^c_j - (\bar{p}_j^M - \bar{d}^u_j - p_j^L)) = 0,
\end{equation}
\begin{equation}\label{kkt4_rh}
-\lambda = c_j \bar{p}^M_j, j \in N,
\end{equation}
\begin{equation}\label{kkt5_rh}
\lambda \in \partial \hat{C}_{h,j}(\bar{d}^c_j), j \in N,
\end{equation}
\begin{equation}\label{kkt6_rh}
\lambda = \bar{d}^u_j/A_j, j \in N.
\end{equation}
\end{subequations}
\end{proposition}

{
\emph{Proof of Proposition \ref{prop_KKT}:}
The proof follows directly from applying subgradient KKT conditions \cite[Section 23]{rockafellar2015convex} to \eqref{Problem_To_Min_RH}.
 \hfill $\blacksquare$
}

\emph{Proof of Theorem \ref{opt_thrm_hybrid}:}
{To prove Theorem \ref{opt_thrm_hybrid}, we {solve} the continuous optimization problem \eqref{Problem_To_Min_RH} using Proposition \ref{prop_KKT}} and then show that {the equilibria are $\epsilon$-optimal to \eqref{Problem_To_Min_RH} which implies that they are also $\epsilon$-optimal to \eqref{Problem_To_Min_H}.}

{First, note that \eqref{kkt5_rh}, i.e. $\lambda \in \partial \hat{C}_{h,j}(\bar{d}^c_j)$, is equivalent to
\begin{equation}\label{dc_optimal}
\bar{d}^c_j = \begin{cases} \bar{d}_j, \lambda > \omega^0_j, \\
(0, \bar{d}_j), \lambda = \omega^0_j, \\
0,  \lambda < \omega^0_j,
\end{cases}  j \in N,
\end{equation}
since $\omega^0_j = \gamma_j$ from \eqref{Design_cond_eqt_1}.
This demonstrates the importance of the constant $\lambda$, which determines the optimum value of on-off load $j$, when $\lambda \neq \omega^0_j$, from \eqref{dc_optimal}.
Furthermore, $\lambda$ needs to be sufficiently large to ensure generation-demand balance, which is reflected in \eqref{kkt1_rh}--\eqref{kkt6_rh}.
Letting $\omega^*$ be the equilibrium frequency of \eqref{sys5}, which is equal for all buses due to~\eqref{eqbr1}, it follows that when $\lambda = \omega^*$, then \eqref{kkt6_rh} holds. Furthermore, when $\alpha_j = c_j^{-1}$, then \eqref{kkt4_rh} also holds from \eqref{eqbr5}. Moreover, {condition \eqref{kkt1_rh} follows} from the summation of equilibrium equation \eqref{eqbr2} over all $j \in N$. Hence, when $\lambda = \omega^*$, if \eqref{dc_optimal} is feasible, i.e. if $d^c_j \in \{0, \bar{d}_j\}, j \in N$,
 then the optimal cost to \eqref{Problem_To_Min_RH} is equal to that of \eqref{Problem_To_Min_H}.
 Below, we explain when \eqref{dc_optimal} is feasible and quantify the additional cost incurred when not.
}


{We denote the minimum costs of the RH-OSLC and H-OSLC problems by $C^{opt}
$ and $C^{\ast}_{opt}$ respectively.
 It then follows that $C^{opt} \leq C^{\ast}_{opt}$
 since the optimal cost to \eqref{Problem_To_Min_RH} provides a lower bound for the global minimum to \eqref{Problem_To_Min_H} as the former is a  relaxed version of the latter, allowing $d^c$ to take continuous values.
Furthermore, let $C^{\ast}$ be the cost associated with \eqref{Problem_To_Min_H} at some equilibrium point {to \eqref{sys5}.}
It then follows that ${C^{\ast} - C^{\ast}_{opt}}  \leq C^{\ast} - C^{opt}$, since $C^{opt} \leq C^{\ast}_{opt} \leq C^{\ast}$.}

 {Note that  Design Condition \ref{Design_condition} allows to deduce the following properties about the equilibria of \eqref{sys5}.
  First, when $p^c \in  [\underline{p}^c_{\underline{j}}, \overline{p}^c_{\underline{j}}]$
  for some $j \in N$, then
  $\sigma^*_{\underline{i}} = 1, i < j, \sigma^*_{\underline{i}} = 0, i > j$ and $\sigma^*_{\underline{j}} \in \{0,1\}$.
 Hence, if $p^c \in \mathcal{F} :=  \bigcup_{j \in N} [\underline{p}^c_j, \overline{p}^c_j]$, then there exist two possible equilibria for $\sigma^*$.
 Alternatively, if $p^c \in \mathbb{R}/\mathcal{F}$ then $\sigma^*$ is unique.
 Note also that $\sigma^*$ determines $\omega^*$ from \eqref{eqlb_freq}
  and  that $[\underline{p}^c_j, \overline{p}^c_j] \cap [\underline{p}^c_k, \overline{p}^c_k] = \emptyset, j \neq k$, as a result of design condition \eqref{Design_cond_eqt_4}.

 Note also that \eqref{kkt_rh} suggests that the value of $\lambda$ is uniquely determined from $\ell$, reminding that $\ell = -p^c$ from \eqref{sys_power_command}.
In particular, when $p^c \in  [\underline{p}^c_{\underline{j}}, \overline{p}^c_{\underline{j}}]$
  for some $j \in N$  then  $\lambda = \omega^0_j$ and when $\lambda \neq \omega^0_j$ for any $j \in N$ then $p^c \in \mathbb{R}/\mathcal{F}$.
 }
{We split the rest of the proof by considering the following two cases: (a) $\lambda \neq \omega^0_j$ for any $j \in N$, (b) there exists $j \in N$ such that $\lambda = \omega^0_j$.}

{\textbf{Part (a):} When  $\lambda \neq \omega^0_j$ for any $j \in N$, then $p^c \in \mathbb{R}/\mathcal{F}$ from Design condition \ref{Design_condition}.
 Hence, as explained above, $\omega^*$ is unique for given $\ell$.
Furthermore, the solution $(\bar{p}^M, \bar{d}^c, \bar{d}^u)$ to \eqref{Problem_To_Min_RH} satisfies $\bar{d}^c_i \in \{0, \bar{d}_i \}, i \in N$ from \eqref{dc_optimal}.
  This suggests that the solutions to \eqref{Problem_To_Min_RH} and \eqref{Problem_To_Min_H} are identical, since $(\bar{p}^M, \bar{d}^c, \bar{d}^u)$ is a feasible solution to \eqref{Problem_To_Min_H} {and therefore $\lambda = \omega^*$.}
   Hence, the equilibria of \eqref{sys5} are global solutions to \eqref{Problem_To_Min_H}.
 }


{\textbf{Part (b):}
As already explained, when a solution to \eqref{Problem_To_Min_RH} satisfies $\lambda = \omega^0_j$ for some $j \in N$,
there exist up to\footnote{{In particular, for given $\ell$, there exist exactly two equilibrium frequencies when $p^c \in \mathcal{F}$ and one otherwise.}}
 two equilibrium frequency values $\omega^*$ to \eqref{sys5} for given $\ell$.
{Furthermore, in general {it can hold that} $\lambda \neq \omega^*$.}
  }}

 {Now let $\lambda = \omega^0_j$ for some $j \in N$, define $S_j = \{l : \omega^0_l = \omega^0_j\}$ and consider a solution $(\bar{p}^M, \bar{d}^c, \bar{d}^u)$ to \eqref{Problem_To_Min_RH}.
 Then, $\bar{d}^c_i \in \{0, \bar{d}_i \}, i \in N /S_j$  {and $\bar{d}^c_i \in [0,\bar{d}_i], i \in S_j$, as follows directly from \eqref{dc_optimal}.} }}
{Now,  the optimal cost to \eqref{Problem_To_Min_RH}, $C^{opt}$, when $\lambda = \omega^0_j$, is given by
\begin{multline*}
{C^{opt} = \sum\limits_{k\in N} \Big( \frac{c_k}{2}(\bar{p}^M_k)^2 +   \hat{C}_{h,k} (\bar{d}^c_k) + \frac{1}{2A_k}(\bar{d}^u_k)^2  \Big)} \\
= \frac{\mathcal{D}}{2} (\omega^0_j)^2 + {\sum\limits_{k\in N} \hat{C}_{h,k} (\bar{d}^c_k).}
\end{multline*}
Analogously, the cost of \eqref{Problem_To_Min_H} at an equilibrium point to \eqref{sys5}, $C^{\ast}$, satisfies
 \begin{multline*}
 C^{\ast}  = {\sum\limits_{k\in N} \Big( \frac{c_k}{2}(p^{M,*}_k)^2 +   C_{h,k} (d^{c,*}_k) + \frac{1}{2A_k}(d^{u,*}_k)^2  \Big)} \\
= \frac{\mathcal{D}}{2} (\omega^*)^2 + {\sum\limits_{k\in N} C_{h,k} ({d}^{c,*}_k)}.
\end{multline*}
Then, note that {$\bar{d}^c_k = {d}^{c,*}_k, k \in N \setminus S_j$}.}
This follows since when $\lambda = \omega^0_j$, then $p^c \in \bigcup_{i \in S_j} [\underline{p}^c_i, \underline{p}^c_i + \overline{d}_i]$ and hence {${d}^{c,*}_k, k \in N \setminus S_j$} satisfy \eqref{dc_optimal} from   \eqref{Design_cond_eqt_2}--\eqref{Design_cond_eqt_4}.
  {Then,
   consider an equilibrium point to \eqref{sys5} and design condition \ref{Design_condition} and note that  for all $\ell$ such that $\lambda = \omega^0_j$ both possible equilibria {satisfy $\hat{q} = \vect{1}^T_{|N|}({d}^{c,\ast} -  \overline{d}^c) \in (-\max_{k \in S_j} \overline{d}_k, \max_{k \in S_j} \overline{d}_k)$. Furthermore, from \eqref{eqlb_freq}, the equilibrium  frequency values
 satisfy $\omega^* = \omega^0_j -\frac{\hat{q}}{\mathcal{D}}$.
  }
Hence,  it follows that  $\sum\limits_{k\in N} (C_{h,k} ({d}^{c,*}_k) - \hat{C}_{h,k} (\bar{d}^c_k)) = \hat{q}\omega^0_j$ since for all $k \in S_j$, the cost per unit demand is $\omega^0_j$ from \eqref{Design_cond_eqt_1}.
 Hence, since $\omega^* = \omega^0_j -\frac{\hat{q}}{\mathcal{D}}$,
 {the difference} between the cost at equilibrium and the optimal cost satisfies
 \begin{multline}\label{cost_diff}
 {C^{\ast} - C^{opt}
 \hspace{-0.5mm} =  \hspace{-0.5mm}   \frac{\mathcal{D}}{2}((\omega^{0}_j)^2 \hspace{-0.5mm}-\hspace{-0.5mm} 2 \frac{\hat{q}}{\mathcal{D}}\omega^0_j \hspace{-0.5mm}+\hspace{-0.5mm} \frac{\hat{q}}{\mathcal{D}^2} \hspace{-0.5mm}-\hspace{-0.5mm} (\omega^{0}_j)^2) \hspace{-0.5mm}+\hspace{-0.5mm}  \hat{q}\omega^0_j }.
\end{multline}
 { Simplifying \eqref{cost_diff} results to $C^{\ast} - C^{opt} = \frac{\hat{q}^2}{2\mathcal{D}}$.}
     Since  $q \in (-\max_{j \in N} \overline{d}_k, \max_{j \in N} \overline{d}_k)$, it follows that {$C^{\ast} - C^{\ast}_{opt} \leq C^{\ast} - C^{opt} \leq
    {\frac{1}{{2\mathcal{D}}}\max_{k \in N} (\bar{d}_k)^2}$},
    which completes the proof.}
 \hfill $\blacksquare$

\balance

\bibliography{andreas_bib}           

\end{document}